\documentclass[a4paper]{article}
\usepackage{authblk}

\renewcommand\Affilfont{\small}
\title{The Intermediate Logic of Convex Polyhedra}
\author[1]{Sam Adam-Day}
\author[2]{Nick Bezhanishvili}
\author[3]{David Gabelaia}
\author[4]{\authorcr Vincenzo Marra}
\affil[1]{Mathematical Institute,\vspace{-0.4em}\authorcr\Affilfont University of Oxford, United Kingdom}
\affil[2]{Institute for Logic, Language and Computation,\vspace{-0.4em}\authorcr\Affilfont University of Amsterdam, The Netherlands}
\affil[3]{A. Razmadze Mathematical Institute,\vspace{-0.4em}\authorcr\Affilfont I. Javakhishvili Tbilisi State University, Georgia}
\affil[4]{Dipartimento di Matematica ``Federigo Enriques'',\vspace{-0.4em}\authorcr\Affilfont Universit\`a degli Studi di Milano, Italy}
\date{\today}

\RequirePackage{silence}
\usepackage[utf8]{inputenc}
\usepackage[hidelinks]{hyperref}
\usepackage[style=alphabetic,maxbibnames=99]{biblatex}
\usepackage{amsmath}
\usepackage{mathtools}
\usepackage{amsthm}
\usepackage{thmtools}
\usepackage{titling}
\usepackage{csquotes}
\usepackage[inline]{enumitem}
\usepackage{nicefrac}
\usepackage[capitalize,noabbrev]{cleveref}
\usepackage{fonttable}
\usepackage{tikz}
\usepackage{tikz-3dplot}
\usepackage[final]{microtype}
\usepackage{setspace}

\usepackage[charter]{mathdesign}
\usepackage[T1]{fontenc}

\usetikzlibrary{positioning,arrows,graphs,arrows.meta,decorations.pathreplacing,calligraphy,cd}

\newtheorem{theorem}{Theorem}[section]
\newtheorem{lemma}[theorem]{Lemma}

\newtheorem{corollary}[theorem]{Corollary}
\theoremstyle{definition}
\newtheorem{definition}[theorem]{Definition}
\theoremstyle{remark}
\newtheorem{remark}[theorem]{Remark}

\tikzset{
	world/.style={draw,circle,outer sep=0pt,inner sep=0,minimum size=15},
	point/.style={draw=black,fill=black,opacity=1,circle,outer sep=0pt,inner sep=0,minimum size=2},
	facefill/.style={fill=blue!20,fill opacity=0.4}
}
\tikzgraphsset{
	poset/.style={grow up=#1, branch right=#1, empty nodes, simple, nodes=point},
	poset/.default=1cm
}

\makeatletter
\tikzdeclarecoordinatesystem{bc}
{%
  {%
    \pgf@xa=0pt
    \pgf@ya=0pt%
    \pgf@xb=0pt
    \tikz@bc@dolist#1,,%
    \pgfmathparse{1/\the\pgf@xb}%
    \global\pgf@x=\pgfmathresult\pgf@xa%
    \global\pgf@y=\pgfmathresult\pgf@ya%
  }%
}%

\def\tikz@bc@dolist#1,{%
  \def\tikz@temp{#1}%
  \ifx\tikz@temp\pgfutil@empty%
  \else
    \pgf@process{\pgfpointanchor{#1}{center}}%
    \pgfmathparse{1}%
    \advance\pgf@xa by\pgfmathresult\pgf@x%
    \advance\pgf@ya by\pgfmathresult\pgf@y%
    \advance\pgf@xb by\pgfmathresult pt%
    \expandafter\tikz@bc@dolist%
  \fi%
}
\makeatother

\newcommand{\Poly}{\ensuremath{\mathsf{Polyhedra}}}
\newcommand{\PolyCon}{\ensuremath{\mathsf{Convex}}}

\newcommand{\SL}{\ensuremath{\mathbf{SL}}}

\newcommand{\IPC}{\ensuremath{\mathbf{IPC}}}
\newcommand{\CPC}{\ensuremath{\mathbf{CPC}}}

\newcommand{\PL}{\ensuremath{\mathbf{PL}}}
\newcommand{\BD}{\ensuremath{\mathbf{BD}}}

\newcommand{\LF}{\ensuremath{\mathbf{LF}}}

\newcommand{\Flat}{\ensuremath{\mathbf{Flat}}}

\DeclareMathOperator{\depth}{\mathsf{depth}}

\DeclareMathOperator{\height}{\mathsf{height}}

\DeclareMathOperator{\Top}{\mathsf{Top}}

\DeclareMathOperator{\Up}{\mathrm{Up}}

\DeclareMathOperator{\Logic}{\mathrm{Logic}}
\DeclareMathOperator{\Frames}{\mathrm{Frames}}
\DeclareMathOperator{\FramesRoot}{\Frames_\bot}
\DeclareMathOperator{\Spec}{\mathrm{Spec}}

\DeclareMathOperator{\uset}{\ua}
\DeclareMathOperator{\dset}{\da}

\DeclareMathOperator{\Opens}{\mathcal{O}}
\DeclareMathOperator{\Closeds}{\mathcal{C}}
\DeclareMathOperator{\Int}{\mathrm{Int}}
\DeclareMathOperator{\Cl}{\mathrm{Cl}}
\newcommand{\comp}[1]{\ensuremath{#1^\mathsf{C}}}

\DeclareMathOperator{\Conv}{\mathrm{Conv}}
\DeclareMathOperator{\Aff}{\mathrm{Aff}}

\DeclareMathOperator{\Relint}{\mathrm{Relint}}
\DeclareMathOperator{\Dim}{\mathrm{Dim}}

\newcommand{\R}{\ensuremath{\mathbb{R}}}

\newcommand{\Subo}{\ensuremath{\mathrm{Sub}_\mathrm{o}}}

\newcommand{\Sub}{\ensuremath{\mathrm{Sub}}}
\newcommand{\Po}{\ensuremath{\mathrm{P}_\mathrm{o}}}
\newcommand{\Pc}{\ensuremath{\mathrm{P}_\mathrm{c}}}

\newcommand{\Lo}{\ensuremath{\mathcal{L}}}
\newcommand{\C}{\ensuremath{\mathbf{C}}}

\renewcommand{\P}{\ensuremath{\mathcal{P}}}

\newcommand{\NN}{\mathbb N}

\newcommand{\N}{\ensuremath{\mathcal{N}}}

\newcommand{\ostar}{\ensuremath{\mathrm{o}}}

\makeatletter
\newcommand{\@usstar}[1]{{\ua}\left(#1\right)}
\newcommand{\@usnostar}[1]{{\ua}(#1)}
\newcommand{\us}{\@ifstar{\@usstar}{\@usnostar}}
\newcommand{\@dsstar}[1]{{\da}\left(#1\right)}
\newcommand{\@dsnostar}[1]{{\da}(#1)}
\newcommand{\ds}{\@ifstar{\@dsstar}{\@dsnostar}}
\newcommand{\@udsstar}[1]{{\uda}\left(#1\right)}
\newcommand{\@udsnostar}[1]{{\uda}(#1)}
\newcommand{\uds}{\@ifstar{\@udsstar}{\@udsnostar}}
\newcommand{\@Usstar}[1]{{\Ua}\left(#1\right)}
\newcommand{\@Usnostar}[1]{{\Ua}(#1)}
\newcommand{\Us}{\@ifstar{\@Usstar}{\@Usnostar}}
\newcommand{\@Dsstar}[1]{{\Da}\left(#1\right)}
\newcommand{\@Dsnostar}[1]{{\Da}(#1)}
\newcommand{\Ds}{\@ifstar{\@Dsstar}{\@Dsnostar}}
\newcommand{\@Udsstar}[1]{{\Uda}\left(#1\right)}
\newcommand{\@Udsnostar}[1]{{\Uda}(#1)}
\newcommand{\Uds}{\@ifstar{\@Udsstar}{\@Udsnostar}}
\makeatother

\newcommand{\Sig}{\Sigma}
\newcommand{\sig}{\sigma}

\newcommand{\vD}{\vDash}

\newcommand{\nvd}{\nvdash}
\newcommand{\nvD}{\nvDash}

\newcommand{\la}{\leftarrow}
\newcommand{\ra}{\rightarrow}

\newcommand{\Lra}{\Leftrightarrow}

\newcommand{\rsa}{\rightsquigarrow}
\newcommand{\cra}{\circrightarrow}
\newcommand{\ua}{\uparrow}
\newcommand{\da}{\downarrow}
\newcommand{\uda}{\updownarrow}
\newcommand{\Ua}{\Uparrow}
\newcommand{\Da}{\Downarrow}
\newcommand{\Uda}{\Updownarrow}
\newcommand{\es}{\varnothing}
\newcommand{\sse}{\subseteq}

\newcommand{\wh}{\widehat}

\newcommand{\half}{\ensuremath{\nicefrac{1}{2}}}

\newcommand{\id}{\ensuremath{\mathsf{id}}}

\newcommand{\abs}[1]{\mathopen|#1\mathclose|}

\newcommand{\contradiction}{\noindent
	\begin{tikzpicture}[x=0.4ex,y=0.4ex]
		\draw[line width=.15ex] (0,0) -- (1,2) -- (0,2) -- (1,4)
		(0.95,0.32) -- (0,0) -- (-0.32,0.95);
	\end{tikzpicture}\hspace*{0.2em}
}
\newcommand{\circrightarrow}{\mathrel{{\circ}\mkern-3mu{\rightarrow}}}

\renewcommand{\leq}{\leqslant}
\renewcommand{\geq}{\geqslant}

\renewcommand{\preceq}{\preccurlyeq}

	{\newline\vspace{\abovedisplayskip}\hbox to \textwidth\bgroup\hss$\displaystyle}
	{$\hss\egroup\vspace{\belowdisplayskip}}

\newcommand{\FOneFork}{
	\raisebox{-0.8ex}{\resizebox{!}{2.5ex}{
		\begin{tikzpicture}
			\node[world] (z) at (0,0) {};
			\node[world] (a) at (0,1) {};
			\path[draw] (z) -- (a);
		\end{tikzpicture}
	}}
}
\newcommand{\FPoint}{
	\raisebox{0.2ex}{\resizebox{!}{0.9ex}{
		\begin{tikzpicture}
			\node[world] (z) at (0,0) {};
		\end{tikzpicture}
	}}
}

\newcommand{\FTwoFork}{
	\raisebox{-0.8ex}{\resizebox{!}{2.5ex}{
		\begin{tikzpicture}
			\node[world] (z) at (0,0) {};
			\node[world] (a1) at (-0.7,1) {};
			\node[world] (b1) at (0.7,1) {};
			\path[draw] (z) -- (a1);
			\path[draw] (z) -- (b1);
		\end{tikzpicture}
	}}
}

\newcommand{\FThreeFork}{
	\raisebox{-0.5ex}{\resizebox{!}{2ex}{
		\begin{tikzpicture}
			\node[world] (z) at (0,0) {};
			\node[world] (a1) at (-1,1) {};
			\node[world] (b1) at (0,1) {};
			\node[world] (c1) at (1,1) {};
			\path[draw] (z) -- (a1);
			\path[draw] (z) -- (b1);
			\path[draw] (z) -- (c1);
		\end{tikzpicture}
	}}
}

\newcommand{\FScott}{
	\raisebox{-1.1ex}{\resizebox{!}{3ex}{
		\begin{tikzpicture}
			\node[world] (z) at (0,0) {};
			\node[world] (a1) at (-0.7,1) {};
			\node[world] (a2) at (-0.7,2) {};
			\node[world] (b1) at (0.7,1) {};
			\path[draw] (z) -- (a1) -- (a2);
			\path[draw] (z) -- (b1);
		\end{tikzpicture}
	}}
}

\DeclareFontFamily{U} {MnSymbolC}{}

\DeclareFontShape{U}{MnSymbolC}{m}{n}{
	<-6> MnSymbolC5
	<6-7> MnSymbolC6
	<7-8> MnSymbolC7
	<8-9> MnSymbolC8
	<9-10> MnSymbolC9
	<10-12> MnSymbolC10
	<12-> MnSymbolC12}{}
\DeclareFontShape{U}{MnSymbolC}{b}{n}{
	<-6> MnSymbolC-Bold5
	<6-7> MnSymbolC-Bold6
	<7-8> MnSymbolC-Bold7
	<8-9> MnSymbolC-Bold8
	<9-10> MnSymbolC-Bold9
	<10-12> MnSymbolC-Bold10
	<12-> MnSymbolC-Bold12}{}

\DeclareSymbolFont{MnSyC} {U} {MnSymbolC}{m}{n}
\DeclareMathSymbol{\meddiamond}{\mathbin}{MnSyC}{110}

\begin{document}

	\maketitle

	\abstract{
		We investigate a recent semantics for intermediate (and modal) logics in terms of polyhedra. The main result is a finite axiomatisation of the intermediate logic of the class of all polytopes --- i.e.,  compact convex polyhedra --- denoted $\PL$. This logic is defined in terms of the Jankov-Fine formulas of two simple frames. Soundness of this axiomatisation requires extracting the geometric constraints imposed on  polyhedra by the two formulas, and then using substantial classical results from polyhedral geometry to show that convex polyhedra satisfy those constraints.  To establish completeness of the axiomatisation, we first define the notion of the geometric realisation of a frame into a polyhedron. We then show that any $\PL$ frame is a p-morphic image of one which has a special form: it is a `sawed tree'. Any sawed tree has a geometric realisation into a convex polyhedron, which completes the proof.
	}

	\renewcommand{\thefootnote}{}
	\renewcommand{\thefootnote}{\arabic{footnote}}
	\addtocounter{footnote}{-2}


\section{Introduction}
\label{sec:introduction}

Polyhedral semantics was introduced in \cite{tarski-polyhedra}. The starting point is that the collection of open subpolyhedra\footnote{For the terminology we adopt in polyhedral geometry the reader is referred to Section \ref{sec:preliminaries}.} of a compact polyhedron (of any dimension) forms a Heyting algebra. This then allows for the interpretation of intuitionistic and modal formulas in polyhedra. This semantics is closely related to the well-known topological semantics, as pioneered in \cite{Stone1938, tsao-chen1938, Tarski1939, mckinsey1941, mckinseytarski44, rasiowasikorski1963}. In topological semantics, one takes the Heyting algebra of open sets of a topological space as the basis for the interpretation of formulas. A celebrated result due to Tarski \cite{Tarski1939} shows that this provides a complete semantics for intuitionistic propositional logic ($\IPC$). The paper \cite{tarski-polyhedra} proved an analogous result for polyhedral semantics: the logic of the class of all polyhedra is \IPC. Moreover, this semantics can access the dimension of a polyhedron via the bounded-depth schema, something beyond the capabilities of topological semantics.

Precursors to the work in \cite{tarski-polyhedra} are \cite{AvBB03, vBBG03, vBB07, interpretating-topological-2010}. In \cite{planar-polygons} and \cite{gabelaiatacl} the authors developed the modal logic of the plane $\R^2$ considered as a non-compact polyhedron. The present authors extended the results of \cite{tarski-polyhedra} in \cite{polycompleteness}, where we introduced the notion of polyhedral completeness: a logic $\Lo$ is polyhedrally complete if it is the logic of some class of polyhedra. We developed the `Nerve Criterion', which provides a necessary and sufficient condition for the polyhedral completeness of a logic based on the combinatorial properties of its frames. This criterion was used to provide a wide class of polyhedrally complete logics axiomatised by the Jankov-Fine formulas of `starlike trees'. The first-named author's M.Sc thesis \cite{sam-thesis} investigated the polyhedral semantics defined in \cite{tarski-polyhedra} and is the basis for both \cite{polycompleteness} and the present paper. Recently, this semantics has been applied to the field of model checking. The authors of \cite{model-checking} developed a geometric spatial model checker using polyhedral semantics, introducing the notion of bisimularity for polyhedra along the way.

In the present paper, we investigate convex (compact) polyhedra, also known as polytopes, from a logical perspective. Our main result (\cref{thm:PL logic of CP}) is that the logic of the class of all convex polyhedra is $\PL$, a logic which is axiomatised by the Jankov-Fine formulas of two simple frames: $\FThreeFork$ and $\FScott$. Moreover, we obtain a more fine-grained result by restricting dimension. Letting $\PL_n$ be $\PL$ plus the logic of bounded depth $n$, we see that this is the logic of the class of all convex polyhedra of dimension at most $n$ (\cref{thm:PLn logic of CPn}).

To prove these results, the first step is a development of the logic-polyhedra connection on the level of morphisms. We introduce the notion of a `polyhedral map' from a polyhedron to a Kripke frame, and show that the open polyhedral maps are exactly those which give rise to contravariant  homomorphisms of the Heyting algebras associated with the polyhedron and the frame, respectively. With this, we can define the notion of the geometric realisation of a frame $F$ to be a polyhedron $P$ together with an open surjective polyhedral map $P \to F$. Moreover, we consider PL (for ``piecewise-linear'') homeomorphisms, which is the standard notion of isomorphism in polyhedral geometry. We show that PL homeomorphisms preserve the logics of polyhedra. 

Now, the proof that $\PL$ is the logic of convex polyhedra consists of two parts: soundness and completeness. For the soundness part, we first make use of the standard geometric fact that every $n$-dimensional convex polyhedron is PL homeomorphic to the $n$-simplex: the `simplest' polyhedron of dimension $n$. Given that PL homeomorphisms preserve logic, it suffices to show that $\PL_n$ is valid on the $n$-simplex, for which we give a geometric proof utilising classical results from polyhedral geometry.

The completeness direction splits into three stages. First, using a combinatorial argument, we show that every $\PL_n$ frame is the p-morphic image of a `sawed tree of height $n$'. This is a frame which has the form of a planar tree with a `saw structure' added on top. Once we have a sawed tree, we show how to realise it geometrically as an $n$-dimensional convex polyhedron. This realisation is built recursively on the frame structure, and makes key use of the fact that sawed trees are planar. Finally, we utilise a result due to Zakharyaschev \cite{zakharyaschev93} which entails that $\PL$ is the intersection of each $\PL_n$, and this completes the proof.

\Cref{sec:preliminaries} introduces the background on intermediate logics and polyhedral geometry, and \cref{sec:polyhedral semantics} introduces polyhedral semantics, following \cite{tarski-polyhedra}. While polyhedra can be used to provide a semantics for both intermediate and modal logics, we focus on the former side here.

\section{Preliminaries}
\label{sec:preliminaries}

	The present paper deals with intermediate logics. In this section we remind the reader of two standard semantics for such logics, and survey the definitions and results which will play their part in what follows. We also present the basic notions of polyhedral geometry that we need in the paper.

\subsection{Posets as Kripke frames}

	A \emph{Kripke frame} for intuitionistic logic is simply a poset $(F,\leq)$. The validity relation $\vD$ between frames and formulas is defined in the usual way. Given a class of frames $\C$, its \emph{logic} is:
	\begin{equation*}
		\Logic(\C) \coloneqq \{\phi \text{ a formula } \mid \forall F \in \C \colon F \vD \phi \}
	\end{equation*}
	Conversely, given a logic $\Lo$, define:
	\begin{equation*}
		\Frames(\Lo) \coloneqq \{F\text{ a Kripke frame} \mid F \vD \Lo\}
	\end{equation*}
	A logic $\Lo$ has the \emph{finite model property} (f.m.p.) if it is the logic of a class of finite frames.

	Fix a poset $F$. For any $x \in F$, its \emph{upset} and \emph{downset} are defined, respectively, as follows.
	\begin{gather*}
		\us x \coloneqq \{y \in F \mid y \geq x\} \\
		\ds x \coloneqq \{y \in F \mid y \leq x\}
	\end{gather*}
	For any set $S \sse F$, its \emph{upset} and \emph{downset} are defined, respectively, as follows.
	\begin{gather*}
		\uset U \coloneqq \bigcup_{x \in U} \us x\\
		\dset U \coloneqq \bigcup_{x \in U} \ds x
	\end{gather*}
	A subframe $U \sse F$ is \emph{upwards-closed} if $U = \uset U$. It is \emph{downwards-closed} if $\dset U = U$. The \emph{Alexandrov topology} on $F$ is the set $\Up F$ of its upwards-closed subsets. This constitutes a topology on $F$. In the sequel, we will freely switch between thinking of $F$ as a poset and as a topological space. Note that the closed sets in this topology correspond to downwards-closed sets.

	A \emph{chain} in $F$ is $X \sse F$ which as a subposet is linearly-ordered. The \emph{length} of the chain $X$ is $\abs X$. A chain $X \sse F$ is maximal if there is no chain $Y \sse F$ such that $X \subset Y$ (i.e. such that $X$ is a proper subset of $Y$). The \emph{height} of $F$ is the element of $\NN \cup \{\infty\}$ defined by:
	\begin{equation*}
		\height(F) \coloneqq \sup\{\abs X-1 \mid X \sse F\text{ is a chain}\}
	\end{equation*}
	For any $x \in F$, define its \emph{height} as follows.
	\begin{gather*}
		\height(x) \coloneqq \height(\ds x)
	\end{gather*}

	The poset $F$ is \emph{rooted} if it has a minimum element, which is called the \emph{root}, and is usually denoted by $\bot$. Define:
	\begin{equation*}
		\FramesRoot(\Lo) \coloneqq \{F \in \Frames(\Lo) \mid F\text{ is rooted}\}
	\end{equation*}

	A function $f \colon F \to G$ is a \emph{p-morphism} if for every $x \in F$ we have:
	\begin{equation*}
		f(\us x) = \us{f(x)}
	\end{equation*}
	An \emph{up-reduction} from $F$ to $G$ is a surjective p-morphism $f$ from an upwards-closed set $U \sse F$ to $G$. Write $f \colon F \cra G$.

	\begin{lemma}\label{lem:up-reduction logic containment}
		If there is an up-reduction $F \cra G$ then $\Logic(F) \sse \Logic(G)$. In other words, if $G \nvD \phi$ then $F \nvD \phi$.
	\end{lemma}

	\begin{proof}
		See \cite[Corollary~2.8, p.~30 and Corollary~2.17, p.~32]{chagrovzakharyaschev1997}.
	\end{proof}

	\begin{corollary}\label{cor:logic of frames logic of rooted frames}
		If $\C$ is any collection of frames and $\Lo = \Logic(\C)$, then:
		\begin{equation*}
			\Lo = \Logic(\FramesRoot(\Lo))
		\end{equation*}
	\end{corollary}

	\begin{proof}
		First, $\Lo \sse \Logic(\FramesRoot(\Lo))$. Conversely, suppose $\Lo \nvd \phi$. Then there exists $F \in \C$ such that $F \nvD \phi$, hence there is $x \in F$ such that $x \nvD \phi$ (for some valuation on $F$), meaning that $\us x \nvD \phi$. Now, $\us x$ is upwards-closed in $F$, hence $\id_{\us x}$ is an up-reduction $F \cra \us x$. Then by \cref{lem:up-reduction logic containment}, we get that $\us x \vD \Lo$, so that $\us x \in \FramesRoot(\Lo)$.
	\end{proof}

	Let $\IPC$ be the logic of all finite frames, and let $\BD_n$ be the logic of all finite frames of height at most $n$. 
	
	\begin{lemma}\label{lem:BDn specifies height}
		Let $F$ be a finite frame. Then $F \vD \BD_n$ if and only if $F$ has height at most $n$.
	\end{lemma}

	\begin{proof}
		See \cite[Proposition~2.38]{chagrovzakharyaschev1997}
	\end{proof}

\subsection{Heyting and co-Heyting algebras}

	A \emph{Heyting algebra} is a tuple $(A,\wedge,\vee,\ra,0,1)$ such that $(A,\wedge,\vee,0,1)$ is a bounded lattice and $\ra$, called the \emph{Heyting implication}, satisfies:
	\begin{equation*}
		c \leq a \ra b \quad\Lra\quad c \wedge a \leq b
	\end{equation*}
	The validity relation $\vD$ between Heyting algebras and formulas is defined in the usual way; the $\Logic$ notation is extended appropriately. Topological spaces provide important examples of Heyting algebras: for every topological space $X$, its collection of open sets $\Opens(X)$ forms a Heyting algebra. We recall that for $U,V\in\Opens(X)$ we have \[U\to V=\bigcup\{Z\in\Opens(X)\mid Z\cap U\subseteq V\}=\Int(U^C\cup V),\]
	where $\Int(-)$ denotes the interior operator and $\comp {(-)}$ denotes set-theoretic complement.

	Co-Heyting algebras are the duals of Heyting algebras. Specifically, a \emph{co-Heyting algebra} is a tuple $(C, \wedge, \vee, \la, 0, 1)$ such that $(C,\wedge,\vee,0,1)$ is a bounded lattice, and $\la$, called the \emph{co-Heyting implication}, satisfies:
	\begin{equation*}
		a \la b \leq c \quad\Lra\quad a \leq b \vee c
	\end{equation*}
	For more information on co-Heyting algebras, the reader is referred to \cite[\S1]{mckinseytarski1946} and \cite{rauszer1974}, where they are called `Brouwerian algebras'.

	Any Heyting algebra $A$ may be regarded as a category. Then its dual category $A^\mathrm{op}$ is a co-Heyting algebra. In the case of the Heyting algebra $\Opens(X)$ of open sets in a topological space, such a duality has a concrete realisation: the co-Heyting algebra $\Opens(X)^\mathrm{op}$ is the algebra $\Closeds(X)$ of closed subsets of $X$.

\subsection{Finite Esakia duality}

	The Alexandrov topology allows us to associate to each poset $F$ the Heyting algebra $\Up F$ consisting of its upwards-closed sets. The process forms part of a contravariant equivalence of categories, known as the Esakia Duality. The finite fragment of this duality relates finite posets with finite Heyting algebras.

	The \emph{spectrum} of a Heyting algebra $A$ is defined as follows. 
	\begin{equation*}
		\Spec(A) \coloneqq \{X \sse A \mid X \text{ is a prime filter of }A\text{ as a distributive lattice}\}
	\end{equation*}
	This constitutes a poset under subset inclusion. 

	\begin{theorem}\label{thm:finite Esakia duality}
		The maps $\Up$ and $\Spec$ are the object-level components of a duality between the category of finite Kripke frames with p-morphisms and the category of finite Heyting algebras with homomorphisms.
	\end{theorem}

	\begin{proof}
		For a proof of the full Esakia Duality see \cite[Corollary~3.4.8]{esakia2019}, which is a translation of the original \cite{esakia1985}. The correspondence was first established in \cite{esakia1974}. Further proofs in English can also be found in \cite{celanijansana2014} and \cite[\S 5]{morandi2005}.
		
		For the finite part, see \cite{dejonghtroelstra1966}. Here, we have isomorphisms $A \cong \Up \Spec A$ and $F \cong \Spec \Up F$ for any finite Heyting algebra $A$ and finite poset $F$. The former is part of Brikhoff's Representation Theorem \cite{birkhoff1937}. Both isomorphisms may be found in \cite[pp.~171-172]{daveypriestley1990}.
	\end{proof}

	Importantly, this duality is logic-preserving.

	\begin{lemma}\label{lem:finite esakia duality logic-preserving}
		Let $F$ be a frame and $A$ be a finite Heyting algebra. Then:
		\begin{gather*}
			\Logic(F) = \Logic(\Up F) \\
			\Logic(A) = \Logic(\Spec A)
		\end{gather*}
	\end{lemma}

	\begin{proof}
		For the first equality, see \cite[Corollary~8.5, p.~238]{chagrovzakharyaschev1997}, noting that our Kripke frames are special cases of what are there called `intuitionistic general frames'. The second equality follows from the first and the finite Esakia duality.
	\end{proof}

\subsection{Jankov-Fine formulas as forbidden configurations}

	To every finite rooted frame $Q$, we associate a formula $\chi(Q)$, the \emph{Jankov-Fine} formula of $Q$ (also called its \emph{Jankov-De Jongh formula}). The precise definition of $\chi(Q)$ is somewhat involved, but the exact details of this syntactical form are not relevant for our considerations. What matters to us is its notable semantic property.

	\begin{theorem}\label{thm:Jankov-Fine up-reductions}
		For any frame $F$, we have that $F \vD \chi(Q)$ if and only if $F$ does not up-reduce to $Q$.
	\end{theorem}

	\begin{proof}
		See \cite[\S9.4, p.~310]{chagrovzakharyaschev1997}, for a treatment in which Jankov-Fine formulas are considered as specific instances of more general `canonical formulas'. A more direct proof is found in \cite[\S3.3, p.~56]{bezhanishvili2006}, which gives a complete definition of $\chi(Q)$. See also \cite{bezhbezh2009} for an algebraic version of this result.
	\end{proof}

	Jankov-Fine formulas formalise the intuition of `forbidden configurations'. The formula $\chi(Q)$ `forbids' the configuration $Q$ from its frames.

\subsection{Polyhedra and simplices}

	Every polyhedron considered here lives in some Euclidean space $\R^n$. Take finitely many points $x_0, \ldots, x_d \in \R^n$. An \emph{affine combination} of $x_0, \ldots, x_d$ is a point $r_0 x_0 + \cdots + r_d x_d$, specified by some $r_0, \ldots, r_d \in \R$ such that $r_0 + \cdots + r_d = 1$. Given a set $S \sse \R^n$, its \emph{affine hull} $\Aff S$ is the collection of affine combinations of its elements. A \emph{convex combination} is an affine combination in which additionally each $r_i \geq 0$. Given a set $S \sse \R^n$, its \emph{convex hull} $\Conv S$ is the collection of convex combinations of its elements. A subset $S \sse \R^n$ is \emph{convex} if $\Conv S = S$. A \emph{polytope} is the convex hull of a finite subset of $\R^n$. A \emph{polyhedron} in $\R^n$ is a set which can be expressed as the finite union of polytopes.  
	
\begin{remark}	A  remark on terminology is in order. In our usage of the term `polyhedron' does not imply convexity, and is the standard one in piecewise-linear topology --- c.f.\@ classic textbooks \cite{stallings1967,rourkesanderson1972}) --- with the following additional conventions. A `polyhedron' \textit{tout court}, as defined in PL topology, need not be compact as a subspace of Euclidean space.  Now, it is a standard fact that `compact polyhedra' (in this more general sense) coincide with what we are referring to in this paper as `polyhedra' (see \cite[Theorem~2.2, p.~12]{rourkesanderson1972}). Hence we are effectively using the term `polyhedron' as a shorthand for `compact polyhedron'. Such abbreviated usage is frequent in the literature (see e.g. \cite{maunder1980algebraic}). Finally,  in our terminology, a `convex polyhedron' is the same thing as a `polytope' --- we will  use the former expression from now on.
\end{remark}
	
	A set of points $x_0, \ldots, x_d$ is \emph{affinely independent} if whenever:
	\begin{equation*}
		r_0 x_0 + \cdots + r_d x_d = \mathbf 0 \quad\text{and}\quad r_0 + \cdots + r_d = 0
	\end{equation*}
	we must have that $r_0=\cdots=r_d = 0$. This is equivalent to saying that the vectors
	\begin{equation*}
		x_1 - x_0, \ldots, x_d - x_0
	\end{equation*}
	are linearly independent. A \emph{$d$-simplex} is the convex hull $\sig$ of $d+1$ affinely independent points $x_0, \ldots, x_d$, which we call its \emph{vertices}. Write $\sig = x_0\cdots x_d$; its \emph{dimension} is $\Dim\sig \coloneqq d$.

	\begin{lemma}
		Every simplex determines its vertex set: two simplices coincide if and only if they share the same vertex set.
	\end{lemma}

	\begin{proof}
		See \cite[Proposition 2.3.3, p.~32]{maunder1980algebraic}.
	\end{proof}

	\noindent A \emph{face} of $\sig$ is the convex hull $\tau$ of some non-empty subset of $\{x_0, \ldots, x_d\}$ (note that $\tau$ is then a simplex too). Write $\tau \preceq \sig$, and $\tau \prec \sig$ if $\tau \neq \sig$.

	Since $x_0, \ldots, x_d$ are affinely independent, every point $x \in \sig$ can be expressed uniquely as a convex combination $x = r_0 x_0 + \cdots + r_d x_d$ with $r_0, \ldots, r_d \geq 0$ and $r_0 + \cdots + r_d = 1$. Call the tuple $(r_0, \ldots, r_d)$ the \emph{barycentric coordinates} of $x$ in $\sig$. The \emph{barycentre} $\wh\sig$ of $\sig$ is the special point whose barycentric coordinates are $(\frac{1}{d+1}, \ldots, \frac{1}{d+1})$. The \emph{relative interior} of $\sig$ is defined as follows.
	\begin{equation*}
		\Relint \sig \coloneqq \{r_0 x_0 + \cdots + r_d x_d \in \sig \mid r_0, \ldots, r_d >0\}
	\end{equation*}
	 Then the relative interior of $\sig$ coincides with the topological interior of $\sig$ inside its affine hull. Note that $\Cl\Relint\sig = \sig$, the closure being taken in the ambient space $\R^n$.

\subsection{Triangulations}

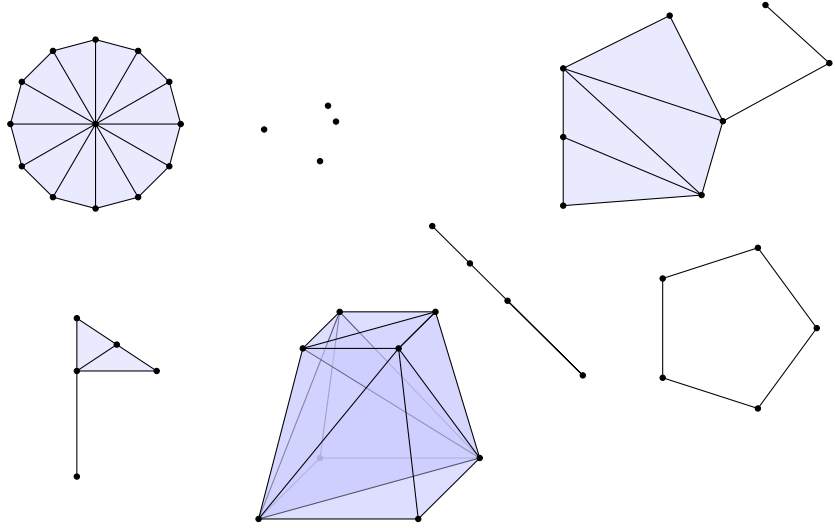
\begin{figure}[ht]
		\begin{equation*}
			\begin{tikzpicture}[scale=0.7]
				\tdplotsetmaincoords{95}{95}
				\begin{scope}[scale=1.6, every node/.style=point]
					\filldraw[facefill]
						(0:1) -- (30:1) -- (60:1) -- (90:1) -- (120:1) -- (150:1) -- (180:1) -- (210:1) -- (240:1) -- (270:1) -- (300:1) -- (330:1) -- cycle;
					\node (z) at (0,0) {};
					\foreach \i in {0,...,11}
					{
						\draw (z) -- ({30*\i}:1) node (\i) {};
					}
				\end{scope}
				\begin{scope}[xshift=250, yshift=30, every node/.style=point]
					\fill[facefill] (0,0) -- (2,1) -- (3,-1) -- (2.6,-2.4) -- (0,-2.6) -- (0,-1.3) -- cycle;
					\draw (0,0) node {} -- (2,1) node {} -- (3,-1) node {} -- (2.6,-2.4) node {} -- (0,-2.6) node {} -- (0,-1.3) node {} -- cycle;
					\draw (3,-1) -- (5,0.1) node {} -- (3.8,1.2) node {};
					\draw (3,-1) -- (0,0) -- (2.6,-2.4) -- (0,-1.3);
				\end{scope}
				\begin{scope}[xshift=120, yshift=-180, scale=3]
					\coordinate (aaa) at (0,0,0);
					\coordinate (baa) at (1,0,0);
					\coordinate (aba) at (0.2,1,0.2);
					\coordinate (bba) at (0.8,1,0.2);
					\coordinate (aab) at (0,0,1);
					\coordinate (bab) at (1,0,1);
					\coordinate (abb) at (0.2,1,0.8);
					\coordinate (bbb) at (0.8,1,0.8);
					\fill[facefill] 	(aaa) -- (aab) -- (bab) -- (baa);
					\fill[facefill] 	(aaa) -- (aab) -- (abb) -- (aba);
					\draw			 	(aaa) -- (aab);
					\fill[facefill] 	(bba) -- (aba) -- (aaa) -- (baa);
					\draw			 	(aba) -- (aaa) -- (baa);
					\node[point] at (0,0,0) {};
					\filldraw[facefill] (aab) -- (aba) -- (baa);
					\fill[facefill]		(baa) -- (abb) -- (bba);
					\fill[facefill]		(baa) -- (abb) -- (aab);
					\draw 				(abb) -- (baa);
					\fill[facefill]		(abb) -- (bbb) -- (baa);
					\fill[facefill] 	(baa) -- (aab) -- (bbb);
					\draw 				(baa) -- (aab);
					\filldraw[facefill] (bbb) -- (bba) -- (baa) -- (bab);
					\filldraw[facefill] (abb) -- (aab) -- (bab) -- (bbb);
					\filldraw[facefill] (aba) -- (abb) -- (bbb) -- (bba) -- cycle;
					\draw 				(baa) -- (bbb) -- (aab);
					\draw				(abb) -- (bba);
					\node[point] at (1,0,0) {};
					\node[point] at (0.2,1,0.2) {};
					\node[point] at (0.8,1,0.2) {};
					\node[point] at (0,0,1) {};
					\node[point] at (1,0,1) {};
					\node[point] at (0.2,1,0.8) {};
					\node[point] at (0.8,1,0.8) {};
				\end{scope}
				\begin{scope}[xshift=120, yshift=-20, scale=1.5]
					\node[point] at (0,0) {};
					\node[point] at (0.2,0.5) {};
					\node[point] at (-0.7,0.4) {};
					\node[point] at (0.1,0.7) {};
				\end{scope}
				\begin{scope}[xshift=180,yshift=-55]
					\draw (0,0) node[point] {}
						-- (-45:1) node[point] {} 
						-- (-45:4) node[point] {} 
						-- (-45:2) node[point] {};
				\end{scope}
				\begin{scope}[xshift=-10,yshift=-190]
					\filldraw[draw=black,facefill]
						(0,2) node[point] (b) {}
							-- ++(0,1) node[point] {}
							-- ++(0.75,-0.5) node[point] (a) {}
							-- ++(0.75,-0.5) node[point] {}
							-- cycle;
					\draw (a) -- (b);
					\draw (0,0) node[point] {} -- (0,2);
				\end{scope}
				\begin{scope}[xshift=340, yshift=-110, scale=1.6]
					\draw (0:1) node[point] {}
						 -- (72:1) node[point] {}
						 -- (144:1) node[point] {}
						 -- (216:1) node[point] {}
						 -- (288:1) node[point] {}
						 -- cycle;
				\end{scope}
			\end{tikzpicture}
		\end{equation*}
		\caption{Triangulations of a collection of polyhedra}
		\label{fig:triangulations of medley of polyhedra}
	\end{figure}

	A \emph{simplicial complex} in $\R^n$ is a finite set $\Sig$ of simplices satisfying the following conditions.
	\begin{enumerate}[label=(\alph*)]
		\item \label{item:downwards-closed; defn:simplicial complex}
		$\Sig$ is $\prec$-downwards-closed: whenever $\sig \in \Sig$ and $\tau \prec \sig$ we have $\tau \in \Sig$.
		\item \label{item:intersection; defn:simplicial complex}
		If $\sig,\tau \in \Sig$, then $\sig \cap \tau$ is either empty or a common face of $\sig$ and $\tau$.
	\end{enumerate}
	The \emph{support} of $\Sig$ is the set $\abs\Sig \coloneqq \bigcup \Sig$. Note that by definition this set is automatically a polyhedron. We say that $\Sig$ is a \emph{triangulation} of the polyhedron $\abs\Sig$. See \cref{fig:triangulations of medley of polyhedra} for some examples of triangulations. 

	Notice that $\Sig$ is a poset under $\prec$, called the \emph{face poset}. A \emph{subcomplex} of $\Sig$ is a subset which is itself a simplicial complex. Note that a subcomplex, as a poset, is precisely a downwards-closed set. Given $\sig \in \Sig$, its \emph{open star} is defined:
	\phantomsection\label{defn:open star}
	\begin{equation*}
		\ostar(\sig) \coloneqq \bigcup \{\Relint \tau \mid \tau \in \Sig \text{ and }\sig \sse \tau\}
	\end{equation*}

		\begin{lemma}\label{lem:open star open}
			The open star $\ostar(\sig)$ of any simplex $\sig$ is open in $\abs\Sig$.
		\end{lemma}

		\begin{proof}
			See \cite[Proposition~2.4.3, p.~43]{maunder1980algebraic}.
		\end{proof}

	\begin{lemma}\label{lem:relints partition support}
		The relative interiors of the simplices in a simplicial complex $\Sig$ partition $\abs\Sig$. That is, for every $x \in \abs\Sig$, there is exactly one $\sig \in \Sig$ such that $x \in \Relint\sig$.
	\end{lemma}

	\begin{proof}
		See \cite[Proposition~2.3.6, p.~33]{maunder1980algebraic}.
	\end{proof}

	In light of \cref{lem:relints partition support}, for any $x \in \abs\Sig$ let us write $\sig^x$ for the unique $\sig \in \Sig$ such that $x \in \Relint\sig$.

	\begin{lemma}\label{lem:simplex relint excludes proper faces}
		Let $\Sig$ be a simplicial complex, take $\tau \in \Sig$ and $x \in \Relint\tau$. Then no proper face $\sig \prec \tau$ contains $x$. This means that $\tau^x = \Relint\tau$ is the inclusion-smallest simplex containing $x$.
	\end{lemma}

	\begin{proof}
		See \cite[Lemma~3.1]{tarski-polyhedra}.
	\end{proof}

	The next result is a basic fact of polyhedral geometry, and is of fundamental importance in its connection with logic. For $\Sig$ a triangulation and $S$ a subspace of the ambient Euclidean space $\R^n$, define:
	\begin{equation*}
		\Sig_S \coloneqq \{\sig \in \Sig \mid \sig \sse S\}
	\end{equation*}
	This, being a downwards-closed subset of $\Sig$, is a subcomplex of $\Sig$.

	\begin{lemma}[Triangulation Lemma]\label{lem:triangulation lemma}
		Any polyhedron admits a triangulation which simultaneously triangulates each of any fixed finite set of subpolyhedra. That is, for a collection of polyhedra $P, Q_1, \ldots, Q_m$ such that each $Q_i \sse P$, there is a triangulation $\Sig$ of $P$ such that $\Sig_{Q_i}$ triangulates $Q_i$ for each $i$.
	\end{lemma}

	\begin{proof}
		See \cite[Theorem~2.11 and Addendum~2.12, p.~16]{rourkesanderson1972}.
	\end{proof}

\subsection{Dimension theory}

	The \emph{dimension} of simplicial complex $\Sig$ is:
	\begin{equation*}
		\Dim\Sig \coloneqq \max\{\Dim\sig \mid \sig \in \Sig\}
	\end{equation*}

	\begin{remark}
		Note that $\Dim\Sig = \height(\Sig)$ as a poset.
	\end{remark}

	\begin{lemma}\label{lem:dimension invariant under triangulation}
		Let $\Sig,\Delta$ be simplicial complexes. If $\abs\Sig = \abs\Delta$ then $\Dim\Sig=\Dim\Delta$.
	\end{lemma}

	\begin{proof}
		See \cite[Proposition~1.6.12, p.~30]{stallings1967}.
	\end{proof}

	\noindent With this in mind, we define the \emph{dimension} $\Dim P$ of a polyhedron $P$ to be the dimension of its triangulations. When $P = \es$, let $\Dim P \coloneqq -1$.

	\begin{lemma}\label{lem:dimension of union}
		$\Dim(P \cup Q) = \max\{\Dim P, \Dim Q\}$.
	\end{lemma}

	\begin{proof}
		By the Triangulation Lemma \ref{lem:triangulation lemma} we can find a triangulation $\Sig$ of $P \cup Q$ such that $\Sig_P$ and $\Sig_Q$ triangulate $P$ and $Q$ respectively. Since $\Sig = \Sig_P \cup \Sig_Q$ and both $\Sig_P$ and $\Sig_Q$ are downwards-closed the result follows.
	\end{proof}

		In the following, it will be necessary to consider the dimensions of sets which are not polyhedra but whose topological closures are. Note that it is possible to define a theory of dimension which applies even more generally \cite{hurewicz-wallmann}, however here we only need to apply it to sets of this form, and the resulting definition is simpler.

		Let $X \sse \R^n$ be such that $\Cl X$ is a polyhedron, where $\Cl X$ denotes the topological closure taken in the ambient space. The \emph{dimension} of $X$ is the dimension of its closure:
		\begin{equation*}
			\Dim X \coloneqq \Dim \Cl X
		\end{equation*}
		\begin{remark}From now on, when we refer to a set $X$ which has dimension, we tacitly assume that its closure is a polyhedron.\end{remark}
		
		Let us  consider the relationship between the dimension operator and the boundary operator. The \emph{boundary} of a set $X$ is $\partial X \coloneqq \Cl^{\Aff} X \setminus \Int^{\Aff} X$, where the closure and interior operations are taken with respect to the affine hull $\Aff X$ (note that $\Cl^{\Aff} X = \Cl X$ in the ambient space, because any affine subspace of $\R^n$ is closed). Then:
		\begin{lemma}\label{lem:dimension of boundary}
			For any set $X$ whose closure is a non-empty polyhedron we have that:
			\begin{equation*}
				\Dim(\partial X) = \Dim(X) - 1
			\end{equation*} 
		\end{lemma}

		\begin{proof}
			See \cite[Corollary~IV.II, p.~46]{hurewicz-wallmann}.
		\end{proof}


\section{Polyhedral semantics}
\label{sec:polyhedral semantics}

With the preliminaries in place, we are in a position to illustrate the link between intuitionistic logic and polyhedra that is the main focus of this paper. Given a polyhedron $P$, let $\Sub P$ denote the collection of its subpolyhedra.

\begin{theorem}\label{thm:subP co-Heyting algebra}
	$\Sub P$ is a co-Heyting algebra, and a subalgebra of $\Closeds(P)$.
\end{theorem}

\begin{proof}
	See \cite[Corollary~3.8]{tarski-polyhedra}. 
\end{proof}

Any subpolyhedron of $P$ is by definition compact, and hence closed. Therefore it is not surprising, once the algebraic nature of $\Sub P$ is established, that it turns out to be a \emph{co-Heyting} algebra. In topology and logic, on the other hand, it is more conventional to work with open sets and \emph{Heyting} algebras. Thus, it is natural at this point to switch to the Heyting algebra dual to $\Sub P$, which has the following concrete realisation.

Given a polyhedron $P$, we will define an \emph{open subpolyhedron} of $P$ as the complement (in $P$) of a  subpolyhedron of $P$; that is, $O\subseteq P$ is an open subpolyhedron of $P$ precisely when the set-theoretic difference $P\setminus O$ is a member of $\Sub P$. 
\begin{remark}Let $P\subseteq \R^n$ be any polyhedron. It is worth pointing out explicitly that while a subpolyhedron of $P$ is a closed (and compact) set both in $P$ and in the ambient space $\R^n$, an open subpolyhedron of $P$ is by definition open in $P$ but may fail to be open in $\R^n$.
\end{remark}
Let us denote by $\Subo P$ the collection of open subpolyhedra in $P$. It is evidently the dual of $\Sub P$, and \cref{thm:subP co-Heyting algebra} yields the following.

\begin{theorem}\label{thm:suboP Heyting algebra}
	$\Subo P$ is a Heyting algebra, and a subalgebra of $\Opens(P)$.
\end{theorem}

The above  provides a sound semantics for intuitionistic logic in terms of polyhedra: for a polyhedron $P$, say that $P \vD \phi$ if and only if $\Subo P \vD \phi$ as a Heyting algebra. One of the features of this polyhedral semantics  is that it is complete for $\IPC$ --- à la Tarski. Moreover, in contrast with topological semantics, polyhedral semantics can detect dimension, via the bounded depth schema. Let $\Poly$ denote the class of all polyhedra, and let $\Poly_n$ denote the subclass consisting of polyhedra of dimension at most $n$, for each $n \in \NN$.

\begin{theorem}\label{thm:IPC logic of P and BDn logic of Pn}
	\begin{enumerate}[label=(\arabic*)]
		\item\label{item:IPC; thm:IPC logic of P and BDn logic of Pn} 
			$\IPC = \Logic(\Poly)$. That is, intuitionistic logic is complete with respect to the class of all polyhedra.
		\item\label{item:BDn; thm:IPC logic of P and BDn logic of Pn} 
			$\BD_n = \Logic(\Poly_n)$, for each $n \in \NN$.
	\end{enumerate}
\end{theorem}

\begin{proof}
	See \cite[Theorem~1.1]{tarski-polyhedra}. The proof works by showing that every finite poset of height $n$ can be `realised geometrically' in an $n$-dimensional polyhedron. The main idea behind this construction is recalled in \cref{ssec:geometric realisation} below.
\end{proof}

	The Triangulation Lemma provides a key piece of information about the polyhedral semantics of Theorem \ref{thm:IPC logic of P and BDn logic of Pn} --- namely, $\Subo P$ is a locally finite Heyting algebra\footnote{An algebraic structure is \emph{locally finite} if every finitely generated substructure is finite.} for any polyhedron $P$. Given any triangulation $\Sig$ of $P$, denote by $\Pc(\Sig)$ the sublattice of $\Closeds(P)$ generated by $\Sig$, and let:
	\begin{equation*}
		\Po(\Sig) \coloneqq \{P \setminus C \mid C \in \Pc(\Sig)\}
	\end{equation*}

	\begin{lemma}\label{lem:Po Sig iso Up Sig}
		$\Po(\Sig)$ is isomorphic as a Heyting algebra to $\Up \Sig$.
	\end{lemma}

	\begin{proof}
		See \cite[Lemma~4.3]{tarski-polyhedra}.
	\end{proof}

	\begin{theorem}\label{thm:Subo P locally-finite}
		Whenever $P \nvD \phi$ there is a triangulation $\Sig$ of $P$ such that $\Po(\Sig) \nvD \phi$. In particular, $\Subo P$ is locally finite.
	\end{theorem}

	\begin{proof}
		See \cite[Corollary~3.7]{tarski-polyhedra}.
	\end{proof}


\section{Logic, polyhedra and morphisms}
\label{sec:polyhedral maps}

In this section we develop assorted functorial aspects of  polyhedral semantics for intermediate logics which are essential ingredients in the main findings of the present paper.

\subsection{Homomorphisms induced by maps of spaces}

We begin with a  result that  requires some preliminary technical definitions.

For $X$ a topological space, by a \emph{lattice basis} for $X$ we mean a sublattice $L$ of the topology $\Opens(X)$ of $X$ that is a basis for that topology. If $L$ is moreover a Heyting subalgebra of the Heyting algebra $\Opens(X)$, we call $L$ a \emph{Heyting basis}.  

If $X$ is a space with a specified Heyting basis $L$ then we define
\[
\Logic(X)\coloneqq \Logic(L),
\]
where in the left-hand side we assume the basis $L$ is understood from context.

For any set $A$, write $\P(A)$ for the complete Boolean algebra of all subsets of $A$. For any function $f\colon A\to B$ between sets, write $f^{-1}\colon\P(B)\to\P(A)$ for the inverse-image function --- given $S\subseteq B$, $f^{-1}[S]\coloneqq\{a \in A \mid f(a)\in S\}$. Then $f^{-1}$ is a homomorphism of Boolean algebras that moreover preserves arbitrary joins and meets.

Now consider spaces $X$ and $Y$ with prescribed lattice bases $L$ and $M$, respectively. A function $f\colon X\to Y$ is \emph{bases-continuous} if $f^{-1}[S]\in L$ for each $S\in M$. Such functions are, of course, continuous. In general, a  function $f\colon X \to Y$ is \emph{open} if $f[U]\in \Opens(Y)$ for each $U \in \Opens(X)$. When $X$ and $Y$ come with prescribed lattice bases $L$ and $M$, let us say that a function $f$ is \emph{bases-open} if $f[U]\in M$ for each $U \in L$. It is clear that such a bases-open  function is open, because the direct-image function $f[-]$ preserves arbitrary unions.

\begin{lemma}\label{lem:maps duality}
	Let  $f \colon X \to Y$ be a function between spaces $X$ and $Y$ with prescribed lattice bases $L$ and $M$, respectively. Write $f^{-1}[-]\colon \P(Y) \to \P(X)$ for the inverse-image function.
	\begin{enumerate}[label=(\arabic*)]
		\item\label{item:a; lem:maps duality} The function $f$ is bases-continuous if and only if $f^{-1}$ descends to a lattice homomorphism $f^* \coloneqq f^{-1}\colon M \to L$. When one of these two equivalent conditions is satisfied, $f$ being surjective implies that $f^*$ is injective.
		\item\label{item:b; lem:maps duality} Assume  further $L$ and $M$ are Heyting bases. Assume the function $f$ is bases-continuous and bases-open. Then $f^{-1}$ descends to a  homomorphism of Heyting algebras $f^*\colon  M \to L$. Moreover, if $f$ is injective then $f^*$ is surjective, and if $f$ is a bijection then $f^*$ is an isomorphism.
	\end{enumerate}
\end{lemma}

\begin{proof}
	Since $f^*$ is a homomorphism of Boolean algebras, the first assertion in \ref{item:a; lem:maps duality} follows from the definitions. For the second assertion in \ref{item:a; lem:maps duality}, suppose $f$ is surjective. Pick $U, V \in M$ distinct, and suppose without loss of generality there is $p \in U \setminus V$. Since $f$ is surjective, there is $x \in X$ with $f(x)=p$. Then $x \in f^{-1}[U]$ but $x \not \in f^{-1}[V]$, so $f^{-1}=f^*$ is injective.
	
	As for \ref{item:b; lem:maps duality}, let us first assume that $f$ is bases-continuous and bases-open, and take $U,V \in M$ with the aim of showing that $f^*(U \ra V) = f^*(U) \ra f^*(V)$. For the left-to-right inclusion, using the fact that $M$ is a basis and that $f^*=f^{-1}[-]$ commutes with Boolean operations, write (letting $\comp S$ denote the complement of $S$): 
	\[
	U\ra V=\Int(\comp U \cup V)=\bigcup\{O\in M \mid O\sse \comp U \cup V\}
	\]
	and:
	\[
	f^{-1}[U] \ra f^{-1}[V]=\Int\left(\comp{ f^{-1}[U]} \cup f^{-1}[V]\right)=\Int \left(f^{-1}[\comp U \cup V]\right).
	\]
	
	Since $f^{-1}[-]$ preserves arbitrary unions too, we obtain $f^{-1}[U\ra V]=\bigcup f^{-1}[O]$ for $O \in M$ ranging over subsets of  $\comp U \cup V$. Now $O\sse \comp U \cup V$ entails $f^{-1}[O]\sse f^{-1}[\comp U\cup V]$. Since $f^{-1}[O]$ is open because $f$ is continuous,  by the definition of interior $f^{-1}[O]\sse \Int(f^{-1}[\comp U\cup V])$, which shows $f^{-1}[U\ra V]\sse f^{-1}[U]\ra f^{-1}[V]$.
		
	For the right-to-left inclusion we have the following chain of inclusions.
	\begin{align*}
		f[f^{-1}[U] \ra f^{-1}[V]]
			&= f\left[\Int\left(\comp{f^{-1}[U]} \cup f^{-1}[V]\right)\right] \\
			&\sse \Int\left(f\left[\comp{f^{-1}[U]} \cup f^{-1}[V]\right]\right) \tag{$f$ is open} \\
			&= \Int\left(f\left[f^{-1}[\comp{U} \cup V]\right]\right) \\
			&\sse \Int(\comp U \cup V) \\
			&= U \ra V
	\end{align*}
	Applying $f^{-1}$ to both sides, we get that $f^{-1}[U] \ra f^{-1}[V] \sse f^{-1}[U \ra V]$. Summing up, $f^*(U \ra V) = f^*(U) \ra f^*(V)$.

	Next, assume $f$ is injective. Let $A\in L$, and let us show $A$ has a pre-image along $f^{*}=f^{-1}$. Certainly  $A\subseteq f^{-1}[f[A]]$. Let us prove the converse inclusion. If $f^{-1}[f[A]]$ is empty then the converse inclusion holds; otherwise, pick $x \in f^{-1}[f[A]]$. Then $f(x)\in f[A]$, so there is $a\in A$ with $f(x)=f(a)$. Since $f$ is injective, $x=a\in A$, and thus $f^{-1}[f[A]]\subseteq A$. Hence $A$ has the pre-image $f[A]$ along $f^{-1}$. Since, moreover, $f$ is bases-open, we have $f[A]\in M$, so $f^*$ is indeed surjective.
	
	Finally, if $f$ is a bijection then by \ref{item:a; lem:maps duality} and what we just proved $f^*$ is a bijective isomorphism of Heyting algebras, and hence an isomorphism.
\end{proof}

\begin{lemma}\label{lem:subspace_quotient}
	Let $X$ be a space, let $L\subseteq \Opens(X)$, let $Y\subseteq X$, and set $M\coloneqq\{O\cap Y\mid O\in L\}$.
	\begin{enumerate}
	\item If $L$ is a (lattice) basis for the topology of $X$ then $M$ is a (lattice) basis for the subspace topology of $Y$.
	\item If $Y$ is open and $L$ is a Heyting basis for the topology of $X$ then $M$ is a Heyting basis for the subspace topology of $Y$.
	\end{enumerate}
\end{lemma}

\begin{proof}
	This is a straightforward verification and shall be omitted.
\end{proof}
	
To deploy Lemmas \ref{lem:maps duality} and \ref{lem:subspace_quotient} in our geometric setting we will require the next fact.

\begin{lemma}\label{lem:convex open subpolyhedra basis}
	The (convex) open subpolyhedra of a (convex) polyhedron $P$ form a basis for the topology on $P$. Moreover, for any polyhedron $P$, $\Subo P$ is a Heyting basis of $P$.
\end{lemma}

\begin{proof}
	Assume $P\subseteq \R^n$ is any polyhedron.
	Take any $x \in P$ and let $U$ be an open neighbourhood of $x$ in $P$. Then there is some open ball $B$ in $\R^n$ about $x$ such that $x \in B \cap P \sse U$. An elementary argument in affine geometry produces  a simplex $\sig$ in $\R^n$ such that $x \in \Relint\sig \sse B$. Then (by the Triangulation Lemma~\ref{lem:triangulation lemma}) the set $\Cl(P \setminus \sig)$ is a compact subpolyhedron of $P$. Its complement $P \cap \Relint \sig$ is therefore an open subpolyhedron of $P$. Furthermore,
	\begin{equation*}
		x \in P \cap \Relint \sig \sse U, 
	\end{equation*}
	which shows   $\Subo P$ is a basis. If $P$ is additionally convex, then $P \cap \Relint \sig$ is also convex because $P$ and $\Relint \sig$ are, which shows that the convex open subpolyhedra of a convex polyhedron form a basis.

	The `moreover' statement follows from the fact that  the basis $\Subo P$ is a Heyting subalgebra of $\Opens(P)$ by Theorem \ref{thm:suboP Heyting algebra}.
\end{proof}

\begin{remark}\label{rem:convention_bases}
	From now on, in light of Lemma \ref{lem:convex open subpolyhedra basis}, we always tacitly assume a polyhedron $P$ is equipped with its Heyting basis $\Subo P$. Also, in light of Lemma \ref{lem:subspace_quotient}, if $Q$ is an open polyhedron in $P$ --- that is, a member of $\Subo P$ for some polyhedron $P$ --- we always tacitly assume that $Q$ is equipped with the Heyting basis $\Subo Q\coloneqq\{O\cap Q \mid O \in \Subo P\}$.
\end{remark}

Finally, in the next definition we isolate the specific instance of basis-con\-ti\-nuous map that is crucial to our context.

\begin{definition}\label{d:polyhedralmap} 
	Let $P$ be a polyhedron and $Y$  a space with a lattice basis $M$. 
	\begin{enumerate*}[label=(\roman*)]
		\item A function $f \colon P \to Y$ is a \emph{polyhedral map} if it is bases-continuous with respect to the  bases $\Subo P$ and $M$, respectively.
		\item Further, let $Q$ be an open subpolyhedron of $P$. A function $f \colon Q \to Y$ is again called a \emph{polyhedral map} if the pre-image of any open set in $M$ is in $\Subo Q$ (see Remark \ref{rem:convention_bases}).
		\item In the special case that the co-domain $Y$ of $f$ is a  poset $F$, we always tacitly assume $M$ is the Heyting basis $\Up F$ of all open sets in the Alexandrov topology on $F$.
		\item When we say  a polyhedral map as in the foregoing items is \emph{open} we always mean it is \emph{bases}-open with respect to the indicated bases.
	\end{enumerate*}
\end{definition}

\subsection{Jankov-Fine, for polyhedra}

\Cref{thm:Jankov-Fine up-reductions} shows that Jankov-Fine formulas encode forbidden configurations for frames. The same is true for polyhedra with respect to polyhedral maps, as we now show.

Let $\Sig$ be a simplicial complex and $F$ a poset. Given any function $f \colon \Sig \to F$, define the map $\wh f \colon \abs\Sig \to F$ by:
\begin{equation*}
	\wh f(x) \coloneqq f(\sig^x)
\end{equation*}

\begin{lemma}\label{lem:p-morphism to open polyhedral map}
	When $f \colon \Sig \to F$ is a p-morphism, $\wh f \colon \abs\Sig \to F$ is an open polyhedral map.
\end{lemma}

\begin{proof}
	For any $U \in \Up F$, we have that:
	\begin{equation*}
		\wh f^{-1}[U] = \bigcup \{\Relint\sig \mid \sig \in \Sig\text{ and } \sig \in f^{-1}[U]\}
	\end{equation*}
	Since $f$ is monotonic, $f^{-1}[U]$ is upwards-closed in $\Sig$ and therefore $\wh f^{-1}[U]$ is an open sub-polyhedron of $\abs\Sig$. Now take an open set $W \sse \abs\Sig$, with the aim of showing that $\wh f[W]$ is open. Define:
	\begin{equation*}
		\Sig\#W \coloneqq \{\sig \in \Sig \mid \Relint(\sig) \cap W \neq \es\}
	\end{equation*}
	Then:
	\begin{equation*}
		\wh f[W] = \{f(\sig^x) \mid x \in W\} = f[\Sig\#W]
	\end{equation*}
	If $\sig \in \Sig\#W$ and $\sig \preceq \tau$, then as $\sig \sse \tau = \Cl\Relint\tau$ and $W$ is open, we have $\tau \in \Sig\#W$; i.e. $\Sig\#W$ is upwards-closed. But now, $f$ is open and so $\wh f[W]$ is also upwards-closed.
\end{proof}

\begin{lemma}\label{lem:Jankov-Fine polyhedral maps}
	Let $P$ be a polyhedron and $F$ a finite rooted frame. Then $P\nvD \chi(F)$ if and only if there exists an open subpolyhedron $Q$ of $P$ and a surjective open polyhedral map $f \colon Q\to F$. Moreover, if $P$ is convex, then we can assume without  loss of generality that $Q$ is also convex. 
\end{lemma}

\begin{proof}
	Let $P\nvD \chi(F)$. By \cref{thm:Subo P locally-finite} there is a triangulation $\Sig$ of $P$ such that $\Po(\Sig) \nvD \chi(F)$, which by \cref{lem:Po Sig iso Up Sig} means that $\Sig \nvD \chi(F)$. Hence by \cref{thm:Jankov-Fine up-reductions} there is an up-reduction $h \colon \Sig \cra F$. Note that $h$ is open (with respect to the Alexandrov topologies) by the definition of p-morphism. Let $H$ be the (upwards-closed) domain of $h$. As $F$ is rooted, $H$ can be assumed without  loss of generality to be rooted --- it suffices to take a pre-image $y$ of the root of $F$ and let $H = \us y$. Applying \cref{lem:p-morphism to open polyhedral map} to the identity map $\id \colon \Sig \to \Sig$ we find an open polyhedral map $\wh\id \colon P \to \Sig$. Let $Q$ be the pre-image of $H$ via $\wh\id$. Then $h\circ \wh\id \colon Q \to F$ is a surjective open polyhedral map.
	
	Now assume that $P$ is convex. Let $x$ be any element in the pre-image of the root of $H$, and note that $Q$ is an open neighbourhood of $x$. Hence by \cref{lem:convex open subpolyhedra basis} there is an open convex subpolyhedron $W\subseteq P$ such that $x \in W\subseteq Q$. Since $\wh\id$ is open, $\wh\id[W]$ is an upwards-closed subset of $H$ containing its root, and therefore $H = \wh\id[W]$. We have thus found a convex open subpolyhedron $W$ such that $h\circ \wh\id[W] = F$, as desired. 
	
	For the converse direction, as $F\nvD \chi(F)$ we obtain from \cref{lem:maps duality} that $Q\nvD \chi(F)$. Then  \cref{lem:subspace_quotient} implies that $\Subo Q$ is a quotient of $\Subo P$ (via the map $O\in P\mapsto O\cap Q\in\Subo Q$), and therefore $P\nvD \chi(F)$. 
\end{proof}

\subsection{PL maps}

For any $X\sse\R^m$, $Y \sse \R^n$, a function $X \to Y$ is an \emph{affine map} if it lifts to a map $\R^m\to\R^n$ of the form $x \mapsto Mx + b$, where $M$ is a linear transformation and $b \in \R^n$. Now let $P$  and $Q$ be polyhedra in $\R^m$ and $\R^n$, respectively. A function $f\colon P \to Q$ is \emph{piecewise linear}, or a \emph{PL map} for short, if there are triangulations $\Sigma$ and $\Delta$ of $P$ and $Q$ respectively such that
\begin{enumerate}[label=(\arabic*)]
	\item the function $f$  agrees on each $\sigma\in\Sigma$ with an affine map, and
	\item for each $\sigma\in\Sigma$, $f[\sigma]\in \Delta$.
\end{enumerate}
PL maps as just defined are automatically continuous. 

\begin{remark}\label{rem:PL_graph} 
	There are several characterisations, or equivalent definitions, of PL map; we mention one that we shall use, referring to \cite{rourkesanderson1972} for proofs: a function $f\colon P \to Q$ is PL if and only if it is continuous, and its graph  $\{(x,f(x)) \in \R^{m+n}\mid x \in P\}$ is a polyhedron.
\end{remark}

\begin{remark}\label{rem:PL_is_poly}
	A PL map is a polyhedral map because of the standard fact that the inverse image of a polyhedron under a PL-map is a polyhedron, cf. \cite[Corollary~2.5, p.~13]{rourkesanderson1972}. The converse is not true --- the map $[0,1]\to[0,1]$ given by $x\mapsto x^2$ is a polyhedral map that is not PL.
\end{remark}

A \emph{PL homeomorphism} is a PL map that is a homeomorphism.

\begin{lemma}\label{lem:inversePL}
	The inverse of a PL homeomorphism is a PL homeomorphism.
\end{lemma}

\begin{proof}
	See \cite[p.~6]{rourkesanderson1972}.
\end{proof}

\begin{corollary}\label{cor:PL homeomorphism HA isomorphism}
	A PL homeomorphism $f\colon P\to Q$ between polyhedra and its inverse $g\colon Q\to P$ induce mutually inverse isomorphisms of Heyting algebras $f^*\coloneqq f^{-1}\colon \Subo{Q} \to \Subo{P}$ and $g^*\coloneqq g^{-1}\colon \Subo{P} \to \Subo{Q}$.
\end{corollary}

\begin{proof}
	This is an immediate consequence of \cref{lem:maps duality} together with Lemma \ref{lem:inversePL} and \cref{rem:PL_is_poly}.
\end{proof}

\begin{corollary}\label{cor:PL homeomorphic implies logics same}
	If $P$ and $Q$ are PL homeomorphic then $\Logic(P)=\Logic(Q)$.
\end{corollary}

\subsection{Geometric realisation}
\label{ssec:geometric realisation}

The notion of `geometric realisation' can now be made more precise. Given a polyhedron and a space $Y$ with a Heyting basis $M$, a  \emph{realisation} of $Y$ in a polyhedron $P$ is an open surjective polyhedral map $f\colon P \to Y$. By  \cref{lem:maps duality} the dual map $f^*\colon M\to \Subo P$ is an injective homomorphism of Heyting algebras, and this entails $\Logic(P) \sse \Logic(Y)\coloneqq \Logic(M)$,  which is the key ingredient in the completeness proofs. 

Let us emphasise that our usage of the term `geometric realisation'  is specific to our setting. The map $f\colon P\to Y$ `realises' the Heyting algebra $M$ as a subalgebra of $\Subo P$  by pulling back inverse images along $f^*\coloneqq f^{-1}$. This applies in particular to the special case in which  $Y$ is a finite poset $F$, and $M$ is $\Up F$.   We shall next show how this notion of realisation for finite posets  relates to the  standard one of geometric realisation of a simplicial complex.

Let us see how to produce a geometric realisation for an arbitrary finite poset $F$ of height $n$, following \cite{tarski-polyhedra}. For this, we make use of the following construction coming from combinatorial geometry. The \emph{nerve} of $F$, denoted $\N(F)$ is the poset of all non-empty chains in $F$ ordered by inclusion. The nerve comes equipped with a p-morphism $\max \colon \N(F) \to F$ which sends a chain to its maximum element. Note also that $\height(\N(F)) = \height(F)$.

Using the nerve, we then define the geometric realisation of $F$ via a simplicial complex. Enumerate $F = \{x_1, \ldots, x_m\}$, and let $e_1, \ldots, e_m$ be the standard basis vectors of $\R^m$. The \emph{simplicial complex induced by} $F$ is defined:
\begin{equation*}
	\nabla F \coloneqq \{\Conv\{e_{i_1}, \ldots, e_{i_k}\} \mid \{x_{i_1}, \ldots, x_{i_k}\} \in \N(F)\}
\end{equation*}
Noting that $\nabla F \cong \N(F)$ as posets, the p-morphism $\max \colon \N(F) \to F$ then induces an open surjective polyhedral map $\abs{\nabla F} \to F$. Furthermore, by definition:
\begin{equation*}
	\Dim \abs{\nabla F} = \height(\N(F)) = n
\end{equation*}
In other words, we have an $n$-dimensional geometric realisation of the height-$n$ poset $F$, which is the main component in the proof of \cref{thm:IPC logic of P and BDn logic of Pn}.

\section{The logic of convex polyhedra}
\label{sec:logic of convex polyhedra}

Recall from \cref{sec:preliminaries} that a polyhedron $P$ is \emph{convex} if $\Conv P = P$, in other words, if the segment joining any two points in $P$ lies entirely in $P$. Let $\PolyCon$ be the class of all convex polyhedra. We can now tackle the question: \emph{what is the logic of all convex polyhedra, $\Logic(\PolyCon)$?} The remainder of the paper will be devoted to a proof that $\Logic(\PolyCon) = \PL$, where $\PL$ is axiomatised by the Jankov-Fine formulas of two simple trees as follows.

\begin{equation*}
	\PL = \IPC + \chi(\FThreeFork) + \chi(\FScott)
\end{equation*}

\begin{theorem}\label{thm:PL logic of CP}
	$\PL$ is the logic of all convex polyhedra: $\PL = \Logic(\PolyCon)$.
\end{theorem}

We show this result by first restricting to the bounded dimension and bounded frame-depth situation, and then use the fact that $\PL$ has the finite model property to obtain the full result. Specifically, let $\PolyCon_n$ denote the class of convex polyhedra of dimension at most $n$, and define:
\begin{equation*}
	\PL_n \coloneqq \BD_n + \PL
\end{equation*}
The main job will be to prove the following.

\begin{theorem}\label{thm:PLn logic of CPn}
	$\PL_n = \Logic(\PolyCon_n)$, for each $n$.
\end{theorem}

\noindent This in turn splits into the following two directions, which will be proved in \cref{sec:soundness} and \cref{sec:completeness}, respectively.

\begin{theorem}[Soundness]\label{thm:PLn soundness for CPn}
	$\PL_n$ is valid on every $P \in \PolyCon_n$.
\end{theorem}

\begin{theorem}[Completeness]\label{thm:PLn completeness for CPn}
	If $\PL_n \nvd \phi$ then there is $P \in \PolyCon_n$ such that $P \nvD \phi$.
\end{theorem}

The final ingredient is the following result due to Zakharyaschev.

\begin{lemma}\label{lem:PL fmp}
	$\PL$ has the finite model property.
\end{lemma}

\begin{proof}
	This follows from the more general result \cite[Corollary 0.11, p.~20]{zakharyaschev93}. This result is stated in terms of `canonical formulas', which are a generalisation of Jankov-Fine formulas. Given a frame $Q$ and a set $\mathfrak D$ of antichains in $Q$ (sets of pairwise incompatible elements of $Q$), we can define the \emph{canonical formula} $\beta(Q, \mathfrak D, \bot)$, which satisfies a similar condition to that satisfied by Jankov-Fine formulas. The result states that if an intermediate logic $\Lo$ is axiomatised by a set of canonical formulas $\beta(Q, \mathfrak D, \bot)$ such that in every $A \in \mathfrak D$ there is at least one point not lying below all maximal points in $\uset A$, then $\Lo$ has the finite model property.
	
	Now, given any frame $Q$, the Jankov-Fine formula $\chi(Q)$ is equivalent to $\beta(Q, \mathfrak D^\#, \bot)$, where $\mathfrak D^\#$ is the set of non-singleton antichains in $Q$ \cite[Proposition 9.41 (i), p.~312]{chagrovzakharyaschev1997}. It is then clear to see that $\chi(\FThreeFork)$ and $\chi(\FScott)$ satisfy the requisite conditions, so the result yields that $\PL$ has the finite model property.
\end{proof}

These lemmas then combine to give the ultimate result.

\begin{proof}[Proof of \cref{thm:PL logic of CP}]
	\cref{lem:PL fmp} entails that:
	\begin{equation*}
		\PL = \bigcap_{n \in \NN} \PL_n
	\end{equation*}
	On the other hand, since all our polyhedra have finite dimension:
	\begin{equation*}
		\PolyCon = \bigcup_{n \in \NN} \PolyCon_n
	\end{equation*}
	Therefore:
	\begin{equation*}
		\Logic(\PolyCon) = \bigcap_{n \in \NN} \Logic(\PolyCon_n)
	\end{equation*}
	\cref{thm:PLn logic of CPn} then completes the proof.
\end{proof}

\subsection{The Logic of a single convex polyhedron}\label{ss:opensimplex}

Any two $n$-simplices $\sigma\subseteq\R^d$ and $\tau\subseteq \R^{d'}$ are PL-homeomorphic --- in fact, affinely homeomorphic. Indeed, since affine maps commute with affine combinations, any bijection of the vertex set of $\sigma$ onto the vertex set of $\tau$ lifts to exactly one bijective affine map $\Aff \sigma \to \Aff \tau$.  Let $e_0, \ldots, e_n$ be the standard basis vectors of $\R^{n+1}$. The \emph{standard $n$-simplex} is $\Delta_n \coloneqq \Conv\{e_0, \ldots, e_n\}$. The following is a classical result.

\begin{lemma}\label{lem:convex n-dimensional polyhedra pl homeomorphic to n-simplex}
	Every $n$-dimensional convex polyhedron is PL-homeomorphic to $\Delta_n$.
\end{lemma}

\begin{proof}
	See \cite[Corollary~2.20, p.~21]{rourkesanderson1972}. There it is shown that \emph{$n$-cells} --- which correspond to our $n$-dimensional convex polyhedra --- are \emph{$n$-balls} --- meaning that they are PL-homeomorphic to the $n$-dimensional cube $[0,1]^n$. Since $\Delta_n$ is a convex polyhedron, the result follows.
\end{proof}

Thus, the logic of all convex polyhedra of dimension at most $n$ is just the logic of any given $n$-dimensional such polyhedron, for instance the   $n$-simplex.

\begin{corollary}\label{cor:logic CPn is logic of n-simplex}
	For any $n$-dimensional convex polyhedron $P$,
	$\Logic(\PolyCon_n) = \Logic(\Delta_n)=\Logic(P)$. 
\end{corollary}

\begin{proof}
	This is immediate from \cref{lem:convex n-dimensional polyhedra pl homeomorphic to n-simplex} using \cref{cor:PL homeomorphic implies logics same}.
\end{proof}

Next, given a convex polyhedron $P$, we are interested in determining the logic of its topological interior in $\Aff P$ --- that is, the logic of a convex open polyhedron of dimension $n$. In the special case that $P$ is an $n$-simplex $\sigma$, its topological interior in $\Aff \sigma$ coincides with its relative interior $\Relint \sigma$.

\begin{lemma}\label{lem:map_of_cubes}
	There exists a surjective open polyhedral map $(0,1)^n\to [0,1]^n$.
\end{lemma}

\begin{proof}
	Let us first assume $n=1$. Consider real numbers $a'<x<a<b<y<b'$. We define a function $f\colon [a',b']\to [x,y]$ by prescribing its action on vertices:
	\[
	f(a')=a, f(b')=b, f(x)=x, f(a)=a,f(b)=b,f(y)=y\,,
	\]
	and  by completing the definition of $f$ through affine extension. Then $f$ is a surjective PL map. Its restriction $g$ to $(a',b')$ is  a polyhedral map that is evidently still surjective onto $[x,y]$, and is moreover open. (To verify $f$ is open let $(\alpha,\beta)\subseteq(a',b')$. If $x\leq \alpha$ and $\beta\leq y$ then $f[(\alpha,\beta)]=(\alpha,\beta)$. If $\alpha\leq x$ and $y\leq \beta$ then $f[(\alpha,\beta)]=[x,y]$. If $\alpha\leq x$ and $\beta\leq y$ then $f[(\alpha,\beta)]=[x,\beta)$. Hence $f$ is open.) This shows the existence of a surjective open polyhedral map $g\colon (0,1)\to [0,1]$ that is the restriction to $(0,1)$ of a PL map $[0,1]\to[0,1]$.

	For $n>1$, consider the product of maps $F\coloneqq f\times\cdots \times f\colon [0,1]^n\to [0,1]^n$ and its restriction to $(0,1)^n$, $G\coloneqq g\times \cdots \times g\colon (0,1)^n\to [0,1]^n$. Then $F$ is PL. Indeed, its graph is the $n$-fold product of copies of the graph of $f$, and the latter graph is a polyhedron because $f$ is PL; hence the graph of $F$ is a polyhedron, too, using the standard fact  that a finite product of polyhedra is a polyhedron. Since $F$ is continuous \cite[Proposition 2.3.6 and p.\ 78]{EngelkingRyszard1989Gt}, and its graph is a polyhedron, then  $F$ is PL (\cref{rem:PL_graph}). This  entails  that $G$ is polyhedral: if $O\in \Subo [0,1]^n$, $F^{-1}[O]\in \Subo [0,1]^n$ because $F$ is PL; then $G^{-1}[O]= F^{-1}[O]\cap (0,1)^n \in\Subo (0,1)^n$. Finally, since a finite product of open maps is open \cite[Proposition 2.3.29]{EngelkingRyszard1989Gt}, $G$ is open.
\end{proof}

\begin{lemma}\label{lem:open_polytope_logic}
	Let $P$ be any convex polyhedron, and let $O$ be its topological interior in $\Aff P$. Then $\Logic(P)=\Logic(O)$.
\end{lemma}

\begin{proof}
	Assume $P$ is of dimension $n$. By \cref{lem:convex n-dimensional polyhedra pl homeomorphic to n-simplex} there is a PL-homeomorphism $f\colon P\to \Delta_n$ with inverse $f^{-1}\colon \Delta_n\to P$ which also is PL (\cref{lem:inversePL}). Hence by \cref{cor:PL homeomorphic implies logics same} we have $\Logic(P)=\Logic(\Delta_n)$. By an elementary topological argument, $f$ and $f^{-1}$ descend to mutually inverse homeomorphisms  $g \colon O \to \Relint \Delta_n$ and $g^{-1}\colon \Relint \Delta_n \to O$. These homeomorphisms are polyhedral because $f$ and $f^{-1}$ are PL. Hence, \cref{lem:maps duality} entails $\Logic{O}=\Logic(\Relint \Delta_n)$. Thus it suffices to prove the lemma for $P=\Delta_n$ and $O=\Relint \Delta_n$.

	The inclusion map $\iota\colon\Relint \Delta_n \to \Delta$ is an injective open polyhedral map, so that its dual $\iota^*\colon \Subo \Delta_n\to \Subo \Relint \Delta_n$ is a surjective homomorphism of Heyting algebras by \cref{lem:maps duality}, which entails $\Logic(\Delta_n)\subseteq\Logic(\Relint \Delta_n)$. For the converse inclusion, \cref{lem:map_of_cubes} and \cref{lem:maps duality} entail $\Logic((0,1)^n)\subseteq \Logic([0,1]^n)$. The  argument in the previous paragraph yields $\Logic(\Delta_n)=\Logic([0,1]^n)$ and $\Logic(\Relint \Delta_n)=\Logic((0,1)^n)$, which completes the proof.
\end{proof}

\subsection{The largest logic}

The importance of convex polyhedra is mirrored on the logical side.

\begin{theorem}\label{thm:largest logic}
	\begin{enumerate}[label=(\arabic*)]
		\item\label{item:infty; thm:largest logic}
		$\PL$ is the largest polyhedrally complete logic of height $\infty$.
		\item\label{item:n; thm:largest logic}
		$\PL_n$ is the largest polyhedrally complete logic of height $n$, for each $n \in \NN$.
	\end{enumerate}
\end{theorem}

The starting point to prove the above theorem  is the observations that every $n$-dimensional polyhedron contains a convex polyhedron of that dimension.

\begin{lemma}\label{lem:convex in every poly}
	If $P$ is $n$-dimensional polyhedron and $m \leq n$ then there is $Q$ an $m$-dimensional convex polyhedron with $Q \subseteq P$.
\end{lemma}

\begin{proof}
	Let $\Sig$ be a triangulation of $P$. Since $P$ has dimension $n$, there is a simplex $\sig \in \Sig$ which has height $m$ (when viewing $\Sig$ as a poset). Then $\sig\subseteq P$ is an $m$-simplex, which is by definition convex.
\end{proof}

The remaining part of the proof  rests on the results of   \cref{ss:opensimplex}.

\begin{proof}[Proof of \cref{thm:largest logic}]
	To prove \ref{item:n; thm:largest logic}, let $\Lo$ be a polyhedrally complete logic of height $n$. Then $\Lo = \Logic(\C)$ for some class $\C$ of polyhedra. We claim that $\C$ contains a polyhedron of dimension at least $n$. Indeed, otherwise $\C \sse \Poly_{n-1}$ so that by \cref{thm:IPC logic of P and BDn logic of Pn} we have:
	\begin{equation*}
		\BD_{n-1} = \Logic(\Poly_{n-1}) \sse \Logic(\C) = \Lo
	\end{equation*}
	By \cref{lem:BDn specifies height} this means that $\Lo$ cannot have frames of height $n$, a contradiction. \contradiction

	So take $P \in \C$ of dimension at least $n$. Then by \cref{lem:convex in every poly} there is $Q$ a convex $n$-dimensional polyhedron with $Q \subseteq P$. Let $O$ be the topological interior of $Q$ in $\Aff Q$. The inclusion $O\subseteq P$ is an open injective polyhedral map, so by \cref{lem:maps duality} we have $\Logic(P) \sse\Logic(O)$. But by \cref{lem:open_polytope_logic} we also have $\Logic(O)=\Logic(Q)$, and by \cref{cor:logic CPn is logic of n-simplex} we know $\Logic(Q) = \PL_n$; hence:
\begin{equation*}
		\Lo = \Logic(\C) \sse \Logic(P) \sse \Logic(O)=\Logic(Q) = \Logic(\Delta_n) = \PL_n
	\end{equation*}

	To prove \ref{item:infty; thm:largest logic}, let $\Lo = \Logic(\C)$ be a polyhedrally complete logic of height $\infty$. We can write $\C = \bigcup_{n \in \NN} \C_n$, where $\C_n = \C \cap \Poly_n$. Then:
	\begin{equation*}
		\Lo = \Logic(\C) = \Logic \left(\bigcup_{n \in \NN} \C_n\right) = \bigcap_{n \in \NN} \Logic (\C_n) \sse \bigcap_{n \in \NN} \PL_n = \PL
	\end{equation*}
	where in the penultimate containment we have used \ref{item:n; thm:largest logic}, and for the last equality we have used that $\PL$ has the finite model property.
\end{proof}

\section{Soundness}
\label{sec:soundness}

The first half of the proof of \cref{thm:PLn logic of CPn} involves showing that:
\begin{equation*}
    \PL_n = \BD_n + \chi(\FThreeFork) + \chi(\FScott)
\end{equation*}
is valid on all of $\PolyCon_n$. The validity of the first summand follows from \cref{thm:IPC logic of P and BDn logic of Pn}, while for the other two we provide geometric arguments utilising classical results about polyhedra and dimension theory.

We first need the following lemma which relates open polyhedral maps to the boundary operation.

\begin{lemma}\label{lem:polyhedral map boundary}
    Let $f$ be a surjective open polyhedral map from $P$ onto a poset $F$. Whenever $x < y$ in $F$ we have  $f^{-1}[x] \sse \partial f^{-1}[y]$.
\end{lemma}

\begin{proof}
    Since $f$ is open and continuous we have:
    \begin{equation*}
        f^{-1}[x]
            \sse f^{-1} [\ds y]
            = f^{-1} [\Cl \{y\}]
            = \Cl f^{-1} [y]
            = \Cl^{\Aff} f^{-1} [y]
    \end{equation*}
    On the other hand $\Int^{\Aff} f^{-1} [y] \sse f^{-1} [y]$ and $f^{-1}[x]$ is disjoint from $f^{-1}[y]$. Hence:
    \begin{equation*}
        f^{-1}[x] \sse \Cl^{\Aff} f^{-1}[y] \setminus \Int^{\Aff} f^{-1}[y] = \partial f^{-1}[y]\qedhere
    \end{equation*}
\end{proof}

Now, the following is a pure dimension-theoretic result, which is essentially the geometric content of the statement that $\PolyCon \vD \chi(\FScott)$.

\begin{lemma}\label{lem:n-2 disconnection}
    Let $X$ be a convex set of dimension\footnote{Recall that whenever we state that a set has a dimension, we implicitly assume that its closure is a polyhedron.} $n$. There is no $Y \sse X$ of dimension $n-2$ or less such that $X \setminus Y$ is disconnected as a subspace of $X$.
\end{lemma}

\begin{proof}
    See \cite[Corollary~IV.1, p.~48]{hurewicz-wallmann}.
\end{proof}

Similarly, the following is essentially the geometric content of $\PolyCon \vD \chi(\FThreeFork)$.

\begin{lemma}\label{lem:n-1 disconnection}
    Let $X$ be a convex set of dimension $n$. There is no $Y \sse X$ of dimension $n-1$ or less such that $X \setminus Y$ can be partitioned into open sets $U$, $V$ and $W$ with $Y \sse \Cl U \cap \Cl V \cap \Cl W$.
\end{lemma}

To prove this we need the following classical result concerning triangulations of convex polyhedra.

\begin{lemma}\label{lem:n-1 simplices in convex polyhedron}
    Let $\Sig$ be a triangulation of a convex $n$-dimensional polyhedron. Then every $(n-1)$-simplex in $\Sig$ is the face of either one or two simplices of $\Sig$. 
\end{lemma}

\begin{proof}
    See \cite[Exercise~II.4, p.~27]{glaser1970geometricalvI}.
\end{proof}

\begin{proof}[Proof of \cref{lem:n-1 disconnection}]
    Assume for a contradiction that $Y$ disconnects $X$ in such a way that $X \setminus Y$ can be partitioned into open sets $U$, $V$ and $W$ with $Y \sse \Cl U \cap \Cl V \cap \Cl W$. By the Triangulation Lemma~\ref{lem:triangulation lemma} take a triangulation $\Sig$ of $\Cl X$ which simultaneously triangulates $\Cl Y$, $\Cl U$, $\Cl V$ and $\Cl W$.

    By \cref{lem:n-2 disconnection} the set $Y$ must have dimension exactly $n-1$. Hence there is an $(n-1)$-simplex $\sig \in \Sig$ such that $\sig \sse \Cl Y$. By \cref{lem:n-1 simplices in convex polyhedron} we have that $\sig$ is the face of either one or two simplices in $\Sig$. Let $\sig$ be the face of $\tau_1$ and $\tau_2$, where we allow that $\tau_1 = \tau_2$. By our choice of $\Sigma$, each $\Relint \tau_i$ is contained in exactly one of $U$, $V$ and $W$. Assume without loss of generality that $\Relint \tau_1 \sse U$. Similarly, assume that either $\Relint \tau_2 \sse U$ or $\Relint \tau_2 \sse V$.

    Now consider the open star of $\sig$:
    \begin{equation*}
        \ostar(\sig) = \Relint \sig \cup \Relint \tau_1 \cup \Relint \tau_2
    \end{equation*}
    By \cref{lem:open star open} this is open in $X$. Since $\ostar(\sig) \cap Y \neq \es$ and $Y \sse \Cl W$ we have that $\ostar(\sig) \cap W \neq \es$. But this is impossible since $\{Y, U, V, W\}$ forms a partition of $X$ and we have $\Relint \sig \sse Y$ and $\Relint \tau_1, \Relint \tau_2 \sse \Cl U \cup \Cl V$. \contradiction
\end{proof}

With all the pieces in place, we are now in a position to prove the desired soundness result.

\begin{proof}[Proof of \cref{thm:PLn soundness for CPn}]
    That $\PolyCon_n \vD \BD_n$ follows by \cref{thm:IPC logic of P and BDn logic of Pn} \ref{item:BDn; thm:IPC logic of P and BDn logic of Pn}.

    To show the validity of $\chi(\FThreeFork)$, suppose for a contradiction that there is a convex polyhedron $P$ such that $P \nvD \chi(\FThreeFork)$. Then by \cref{lem:Jankov-Fine polyhedral maps} there is a convex open subpolyhedron $Q$ of $P$ and a surjective open polyhedral map $f \colon Q \to \FThreeFork$. By \cref{lem:polyhedral map boundary} this partitions $Q$ into subsets $X, U, V, W$ such that $U$, $V$ and $W$ are open subpolyhedra of $P$ and:
    \begin{equation*}
        X \sse \partial U, \quad
        X \sse \partial V, \quad
        X \sse \partial W
    \end{equation*}
    By \cref{lem:dimension of boundary} we have that $\Dim X \leq \Dim Q - 1$ but $Q \setminus X = U \cup V \cup W$ is disconnected with at least three connected components. This contradicts \cref{lem:n-1 disconnection}. \contradiction

    As for the validity of $\chi(\FScott)$, suppose again for a contradiction that there is a convex polyhedron $P$ such that $P \nvD \chi(\FScott)$.  By \cref{lem:Jankov-Fine polyhedral maps} there is a convex open subpolyhedron $Q$ of $P$ and a surjective open polyhedral map $f \colon Q \to \FScott$.
    Then by \cref{lem:polyhedral map boundary} this partitions $Q$ into subsets $X, U_1, U_2, V_1$ such that $U_1$ and $V_1$ are open subpolyhedra of $P$ and:
    \begin{equation*}
        X \sse \partial U_1, \quad
        U_1 \sse \partial U_2, \quad
        X \sse \partial V_1
    \end{equation*}
    By \cref{lem:dimension of boundary} we have that $\Dim X \leq \Dim Q - 2$ but $Q \setminus X = (U_1 \cup U_2) \cup V_1$ is disconnected. This contradicts \cref{lem:n-2 disconnection}. \contradiction
\end{proof}


\section{Completeness}
\label{sec:completeness}

	The proof that $\PL_n$ is complete with respect to the class of convex polyhedra of dimension at most $n$ consists of two main parts. In the first part, we show that $\PL_n$ can be expressed as the logic of a set of reasonably regular finite frames --- called \emph{sawed trees}. For the second part, we show that any such sawed tree of height $n$ can be realised geometrically as an $n$-dimensional convex polyhedron --- in other words, given a sawed tree $F$, we construct an open polyhedral map from a convex polyhedron onto $F$. This map is constructed using a more elaborate version of the method used to provide a geometric realisation for an arbitrary finite poset in \cref{ssec:geometric realisation}.

\subsection{The meaning of \texorpdfstring{$\PL_n$}{PLn} on frames}

	First of all, it will be convenient to spell out what it means, structurally, for a frame to satisfy $\PL_n$. For this we introduce some additional terminology and notation.

	For any poset $F$ and $x \in F$, the \emph{strict upset} and \emph{strict downset} are defined, respectively, as follows.
	\begin{gather*}
		\Us x \coloneqq \{y \in F \mid y > x\} \\
		\Ds x \coloneqq \{y \in F \mid y < x\}
	\end{gather*}
	The \emph{depth} of $x$ is defined:
	\begin{equation*}
		\depth(x) \coloneqq \height(\us x)
	\end{equation*}
	A \emph{top element} of $F$ is $t \in F$ such that $\depth(t)=0$. The set of top elements in $F$ is denoted by $\Top(F)$.

	A \emph{path} in $F$ is a sequence $p=x_0\cdots x_k$ of elements of $F$ such that for each $i$ we have $x_i < x_{i+1}$ or $x_i > x_{i+1}$. Write $p \colon x_0 \rsa x_k$. The poset $F$ is \emph{path-connected} if between any two points there is a path.

	\begin{lemma}\label{lem:finite poset path-connected iff connected}
		When $F$ is finite, it is path-connected if and only if it is connected as a topological space.
	\end{lemma}

	\begin{proof}
		See \cite[Lemma~3.4]{Bezhanishviligabelaia2011}.
	\end{proof}

	A \emph{connected component} of $F$ is a subframe $U \sse F$ which is connected as a topological subspace and is such that there is no connected $V$ with $U \subset V$.

	\begin{lemma}\label{lem:properties of connectedness}
		\begin{enumerate}[label=(\arabic*)]
			\item The connected components partition $F$.
			\item Connected components are upwards- and downwards-closed.
		\end{enumerate}
	\end{lemma}

	\begin{proof}
		The first is a standard fact in topology, while the second follows straightforwardly from the fact that by \cref{lem:finite poset path-connected iff connected} the connected components are exactly the equivalence classes under the relation `there is a path from $x$ to $y$'.
	\end{proof}

	Finally, for any $x,y \in F$, say that $x$ is an \emph{immediate predecessor} of $y$ and that $y$ is an \emph{immediate successor} of $x$ if $x < y$ and there is no $z \in F$ such that $x < z < y$.

	We can now describe the structural meaning of $\PL_n$ on frames.

	\begin{lemma}\label{lem:meaning of PLn on posets}
		Let $F$ be a poset. Then $F \vD \PL_n$ if and only if the following are satisfied.
		\begin{enumerate}[label=(\roman*)]
			\item\label{item:height; lem:meaning of PLn on posets}
				$F$ has height at most $n$.
			\item\label{item:depth 1; lem:meaning of PLn on posets}
				Whenever $\depth(x) = 1$, we have $\abs{\Us x} \leq 2$.
			\item\label{item:depth gt1; lem:meaning of PLn on posets}
				Whenever $\depth(x) > 1$, the set $\Us x$ is connected.
		\end{enumerate}
	\end{lemma}

	\begin{proof}
		This follows from the definition of $\PL_n$, using the following facts for finite frames $F$.
		\begin{enumerate}[label=(\roman*)]
			\item $F \vD \BD_n$ if and only if $F$ has height at most $n$.
			\item There is an up-reduction $F \cra \FThreeFork$ if and only if there is $x \in F$ such that $\Us x$ has at least three components.
			\item There is an up-reduction $F \cra \FScott$ if and only if there is $x \in F$ such that $\Us x$ has at least two components, with at least one of which having height greater than $0$.\qedhere
		\end{enumerate}
	\end{proof}

	$\PL_n$-frames also satisfy the following specific connectedness property, which will come in handy in the arguments below.

	\begin{lemma}\label{lem:PLn frames specific connectedness property}
		Let $F$ be a finite rooted frame with $\height(F) > 1$, such that $F \vD \PL_n$. Take $s,t \in \Top(F)$. There is a path $p = a_0 \cdots a_m$ from $s$ to $t$ in $\Us{\bot}$ with the property that for each $i$:
		\begin{enumerate}[label=(\Roman*)]
			\item\label{item:a; lem:PLn frames specific connectedness property} 
				$\Us{a_i} = \es$ when $i$ is even, and
			\item\label{item:b; lem:PLn frames specific connectedness property} 
				$\Us{a_i} = \{a_{i-1},a_{i+1}\}$ when $i$ is odd.
		\end{enumerate}
	\end{lemma}

	\begin{proof}
		Since $\height(F) > 1$ we have that $\depth(\bot) > 1$. Hence by \cref{lem:meaning of PLn on posets}, there is a path $p = a_0 \cdots a_m$ from $s$ to $t$ in $\Us{\bot}$. We may assume that:
		\begin{enumerate}[label=(\Alph*)]
			\item\label{item:immediate; proof:PLn frames specific connectedness property} 
				$a_{i+1}$ is either an immediate successor or an immediate predecessor of $a_i$, for each $i$,
			\item\label{item:height-maximal; proof:PLn frames specific connectedness property} 
				$p$ is `height-maximal': if $i < j < k$ and $a_j < a_i, a_k$, then there is no path $a_i \rsa a_k$ in $\Us{a_j}$, and
			\item\label{item:no repeats; proof:PLn frames specific connectedness property} 
				$p$ has no repeats.
		\end{enumerate}
		Indeed, \ref{item:height-maximal; proof:PLn frames specific connectedness property} can be secured by iteratively replacing each offending $a_j$ with the path $a_i \rsa a_k$ in $\Us{a_j}$. Then \ref{item:no repeats; proof:PLn frames specific connectedness property} can be secured by removing all cycles, a process which preserves \ref{item:height-maximal; proof:PLn frames specific connectedness property}.

		We claim that such a $p$ also satisfies \ref{item:a; lem:PLn frames specific connectedness property} and \ref{item:b; lem:PLn frames specific connectedness property}, which we prove by induction. The base $i=0$ is immediate since $a_0=s$ is a top node. So assume that $i>0$. The first case is when $i$ is odd. By induction hypothesis $\Us{a_{i-1}} = \es$; in other words $a_{i-1}$ is a top node. Hence by \ref{item:immediate; proof:PLn frames specific connectedness property}, $a_i$ is an immediate predecessor of $a_{i-1}$. This means that $\{a_{i-1}\}$ is a connected component in $\Us{a_i}$, and hence by \cref{lem:meaning of PLn on posets} \ref{item:depth 1; lem:meaning of PLn on posets} and \ref{item:depth gt1; lem:meaning of PLn on posets}, we must have $\abs{\Us{a_i}} \leq 2$. Note further that by \ref{item:height-maximal; proof:PLn frames specific connectedness property}, $a_{i+1} \neq a_{i-1}$. Therefore, the task is to show that $a_{i+1} \in \Us{a_i}$. Let us suppose for a contradiction that this is not the case; i.e. $a_{i+1} < a_i$. Since $t$ is a top node, there must be $j \geq i+1$ with $a_j \leq a_{i+1}$ such that $a_{j+1} > a_j$ (in other words, the path can not keep going downwards after $a_{i+1}$). Clearly $\depth(a_j)>1$, hence by \cref{lem:meaning of PLn on posets} \ref{item:depth gt1; lem:meaning of PLn on posets} there must be a path $a_i \rsa a_{j+1}$ in $\Us{a_j}$, which contradicts property \ref{item:height-maximal; proof:PLn frames specific connectedness property}. \contradiction Thus $a_{i+1} \in \Us{a_i}$ as required. The second case when $i$ is even follows immediately from property \ref{item:immediate; proof:PLn frames specific connectedness property} and the induction hypothesis.
	\end{proof}

\subsection{Sawed trees}

\begin{figure}[ht]
		\begin{equation*}
			\begin{tikzpicture}[scale=0.8]
				\begin{scope}
					\graph[poset=0.8] {
						a1[x=4cm] -- { 
							b1[x=1cm] -- { 
								c1[x=0.5cm] -- {
									d1 -- {e2[x=0.5cm]},
									d2 -- {e2, e3[x=0.5cm]}
								}, 
								c2 -- d3 -- {e3, e4[x=0.5cm]}
							},
							b2[x=1cm] -- {
								c3 -- d4 -- {e4, e5[x=0.5cm]},
								c4 -- d5 -- {e5, e6[x=0.5cm]},
								c5 -- d6 -- {e6, e7[x=0.5cm]},
							},
							b3 -- c6 -- d7 -- {e7, e8[x=0.5cm]},
							b4 -- c7[x=0.5cm] -- {
								d8 -- {e8, e9[x=0.5cm]},
								d9 -- {e9}
							}
						};
					};
					\draw[dotted,rounded corners=15] (-0.3,4.3) rectangle (8.3,2.7);
					\draw[dashed,rounded corners=20] (4,-0.4) -- (7.3,0.7) -- (8.6,3.35) -- (-0.6,3.35) --  (0.7,0.7) -- cycle;
					\node at (9.3,3.5) {Saw};
					\node at (8.3,1) {Tree};
				\end{scope}
			\end{tikzpicture}
		\end{equation*}
		\caption{An example sawed tree}
		\label{fig:sawed tree example}
	\end{figure}
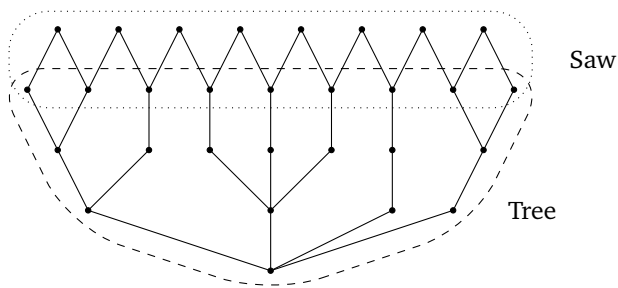

	Let $T$ be a finite tree in which every top element has the same height. A linear ordering $\prec$ on $\Top(T)$ (or equivalently an enumeration $t_1, \ldots, t_k$ of $\Top(T)$) is a \emph{plane ordering} if for every $x \in T$ we have that $\us x \cap \Top(T)$ is an interval with respect to $\prec$. When $\height(T)>0$, the \emph{sawed tree} based on $(T,\prec)$ consists of $T$ plus new elements $s_1, \ldots, s_{k-1}$ with relations, for each $i$:
	\begin{equation*}
		t_i,t_{i+1} < s_i
	\end{equation*}
	See \cref{fig:sawed tree example} for an example of a sawed tree.

	The planarity condition on $\prec$ ensures that the Hasse diagram of the resulting sawed tree can be drawn in the plane with no overlapping lines. Formally, let $G$ be a poset and $d \colon G \to \R^2$ be an injection, such that $d = (d_1,d_2)$. Draw an edge $\mathit{xy}$ between $d(x)$ and $d(y)$ whenever $y$ is an immediate successor of $x$. Then $d$ is a \emph{plane drawing} of $G$ if the following conditions hold.
	\begin{enumerate}[label=(\alph*)]
		\item Whenever $x < y$ we have $d_2(x) < d_2(y)$.
		\item Two distinct edges $x_1y_1$ and $x_2y_2$ only ever intersect at their end-points.
	\end{enumerate}
	The notion of a planar poset has been studied somewhat in the literature (see \cite[\S6.8, p.~101]{brandstadtetal1999} for a short survey), but we will not  use any external results here.

	\begin{lemma}\label{lem:trees are planar}
		Let $\prec$ be a plane ordering on $T$. Then $T$ has plane drawing $d$ with the following properties. 
		\begin{enumerate}[label=(\roman*)]
			\item The top nodes in the drawing are ordered left-to-right as per $\prec$.
			\item $d_2(x) = \height(x)$ for every $x \in T$.
		\end{enumerate}
	\end{lemma}

	\begin{proof}
		C.f. \cite[p.~294]{stanley1997}. We proceed by induction on $n = \height(T)$. The base case $n=0$ is immediate, so assume that $n>0$. Enumerate the immediate successors of $\bot$ in $T$ as $\{x_1, \ldots, x_k\}$, according to $\prec$. That is, for each $i,j \leq k$ with $i<j$ ensure that:
		\begin{equation*}
			\forall t_i \in \us{x_i} \cap \Top(T) \colon \forall t_j \in \us{x_j} \cap \Top(T) \colon t_i \prec t_j
		\end{equation*}
		This is possible since $\us{x} \cap \Top(T)$ is an interval for each $x$. By induction hypothesis, for each $i \leq k$ there is a plane drawing $d^i$ of $\us{x_i}$ satisfying the conditions. We can then form a plane drawing $d$ of $T$ by shifting the drawings $d_1, \ldots, d_k$ up by one, lining them up side by side, then letting $d(\bot) \coloneqq (0,0)$. It is clear that $d$ then also satisfies the required conditions.
	\end{proof}

	\begin{corollary}\label{cor:sawed trees are planar}
		Every sawed tree $F$ admits a plane drawing $d$ with the property that $d_2(x) = \height(x)$ for every $x \in F$.
	\end{corollary}

	\begin{proof}
		Let $F$ be based on $(T,\prec)$, and let $s_1, \ldots, s_{k-1}$ be the top elements. By \cref{lem:trees are planar}, there is a plane drawing $d'$ of $T$ satisfying the property. Extend $d'$ to a drawing $d$ of $F$ by letting $d(s_i) \coloneqq (i,\height(F))$.
	\end{proof}

	The reason for considering sawed trees is that they provide a complete class of frames for $\PL$ which is relatively easy to work with.

	\begin{lemma}\label{lem:sawed trees satisfy PLn}
		Let $F$ be a sawed tree of height $n$. Then $F \vD \PL_n$.
	\end{lemma}

	\begin{proof}
		Let $F$ be based on $(T,\prec)$. Let us verify the conditions of \cref{lem:meaning of PLn on posets}. Conditions \ref{item:height; lem:meaning of PLn on posets} and \ref{item:depth 1; lem:meaning of PLn on posets} are immediate. As for \ref{item:depth gt1; lem:meaning of PLn on posets}, take $x \in F$ with $\depth(x)>1$. By construction, $x \in T$. Since $\prec$ is a plane ordering, we have that $\us x \cap \Top(T)$ is an interval with respect to $\prec$. Therefore, the top two layers of $\Us x$ are connected by the saw structure.
	\end{proof}

	\begin{lemma}\label{lem:frames of PL p-morphic images of sawed trees}
		Every rooted frame $F$ of $\PL$ of height $n$ is the p-morphic image of a sawed tree of height $n$, for every $n \geq 2$.
	\end{lemma}

	\begin{proof}
		We prove this by induction on $n$. For the base case $n=2$, note that $F$ consists of the root $\bot$ together with a number of nodes of depths $0$ and $1$. By gluing together paths obtained from \cref{lem:PLn frames specific connectedness property}, we can find a path $p = a_0 \cdots a_m$ satisfying \ref{item:a; lem:PLn frames specific connectedness property} and \ref{item:b; lem:PLn frames specific connectedness property} of that lemma which visits every top node. We would like to extend $p$ so that it visits \emph{every non-root} node. To do this, take $x \in F$ of depth $1$. By \cref{lem:meaning of PLn on posets} \ref{item:depth 1; lem:meaning of PLn on posets}, $\Us{x}=\{s,t\}$ with $s,t$ top nodes and possibly $s=t$. By inserting the sequence $xtxs$ in $p$ after an occurrence of $s$, we obtain a path satisfying \ref{item:a; lem:PLn frames specific connectedness property} and \ref{item:b; lem:PLn frames specific connectedness property}, which also visits $x$.

		Therefore, we may assume that our path $p$ visits every non-root node. Now, construct the sawed tree $F'$ by taking $\bot$ together with new elements:
		\begin{equation*}
			w_{-1},w_0, \ldots, w_m, w_{m-1}
		\end{equation*}
		with relations as in \cref{fig:F prime n2; proof:frames of PL p-morphic images of sawed trees}.

		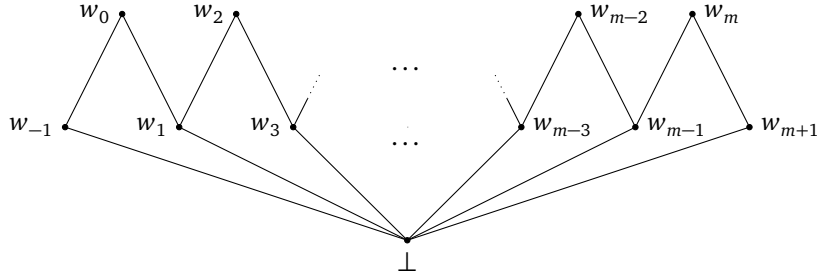
\begin{figure}[ht]
			\begin{equation*}
				\begin{tikzpicture}[yscale=1.5,xscale=1.5]
					\graph[poset=1.5cm] {
						a[x=3cm, label=below:$\bot$] -- {
							b0[label=left:$w_{-1}$] -- c0[x=0.5cm, label=left:$w_0$],
							b1[label=left:$w_1$] -- {c0, c1[x=0.5cm, label=left:$w_2$]},
							b2[label=left:$w_3$] -- c1,
							-!- e0[minimum size=0],
							b3[label=right:$w_{m-3}$] -- c2[x=0.5cm, label=right:$w_{m-2}$],
							b4[label=right:$w_{m-1}$] -- {c2, c3[x=0.5cm, label=right:$w_{m}$]},
							b5[label=right:$w_{m+1}$] -- c3
						}
					};
					\draw (b2) -- ++(0.125,0.25) edge[dotted] +(0.125,0.25);
					\draw (b3) -- ++(-0.125,0.25) edge[dotted] +(-0.125,0.25);
					\node[above=1 of a] {$\cdots$};
					\node[above=2 of a] {$\cdots$};
				\end{tikzpicture}
			\end{equation*}
			\caption{The relations in $F'$ when $n=2$}
			\label{fig:F prime n2; proof:frames of PL p-morphic images of sawed trees}
		\end{figure}

		Then define the surjective map $f \colon F' \to F$ by:
		\begin{gather*}
			\bot \mapsto \bot, \\
			w_{-1} \mapsto a_0, \\ 
			w_{m+1} \mapsto a_m, \\
			w_i \mapsto a_i \qquad \forall i \in \{0, \ldots, m\}
		\end{gather*}
		That $f$ is a p-morphism amounts to the fact that $p$ satisfies properties \ref{item:a; lem:PLn frames specific connectedness property} and \ref{item:b; lem:PLn frames specific connectedness property} of \cref{lem:PLn frames specific connectedness property}.

		For the induction step, assume that $n > 2$. Let $z_1, \ldots, z_k$ be the immediate successors of $\bot$ in $F$. By induction hypothesis, for each $i$ there is a sawed tree $G_i$ and a p-morphism $g_i \colon G_i \to \us{z_i}$. Let the sawed tree $G_i$ be based on $(S_i,\prec_i)$, and let $u_i, v_i \in \Top(S_i)$ be the least and greatest elements according to $\prec_i$, respectively. Since $\abs{\us{u_i}}, \abs{\us{v_i}} = 2$, we must have: 
		\begin{equation*}
			\abs{\us{g_i(u_i)}}, \abs{\us{g_i(v_i)}} \leq 2
		\end{equation*}
		Let $s_i \in \us{g_i(u_i)}$ and $t_i \in \us{g_i(v_i)}$ be the greatest elements. Now, by \cref{lem:PLn frames specific connectedness property}, for each $i \leq k-1$ there is a path $p_i \colon t_i \rsa s_{i+1}$ satisfying properties \ref{item:a; lem:PLn frames specific connectedness property} and \ref{item:b; lem:PLn frames specific connectedness property}; write $p_i = a_{i,0} \cdots a_{i,m_i}$.

		We will form our new sawed tree by laying the sawed trees $G_1, \ldots, G_k$ in a line and `gluing' them usings the paths $p_1, \ldots, p_{k-1}$ together with some `rope ladders' beneath. In detail, form $F'$ by taking the following ingredients and combining them as in \cref{fig:F prime construction inductive step; proof:frames of PL p-morphic images of sawed trees}.
		\begin{itemize}
			\item Each sawed tree $G_i$.
			\item For each $i \leq k$, new elements $w_{i,0} \cdots w_{i,k_i}$ corresponding to $a_{i,0} \cdots a_{i,k_i}$.
			\item A chain of length $n-2$ (a rope ladder) to hang below each $w_{i,j}$, with $j$ odd.
		\end{itemize}

		\begin{figure}[ht]
			\begin{equation*}
				\begin{tikzpicture}[xscale=0.9]

					\newcommand{\sawedtreeblock}[1]{
						\begin{scope}[xscale=0.15, yscale=0.4]
							\graph[grow right=0.135cm, empty nodes, nodes={inner sep=0, outer sep=0, minimum size=0}] {
								#1_w_0
									-- #1_a_0[y=1cm]
									-- #1_w_1
									-- #1_a_1[y=1cm]
									-!- #1_e0[minimum size=0]
									-!- #1_ec[x=0.5cm, y=0.25cm, minimum size=0]
									-!- #1_e1[minimum size=0]
									-!- #1_e2[minimum size=0]
									-!- #1_a_{m-2}[y=1cm]
									-- #1_w_{m-1}
									-- #1_a_{m-1}[y=1cm]
									-- #1_w_m
									;
							};
							\node at (#1_ec) {$\cdots$};
							\node at (5.5,-2) {$G_{#1}$};
							\draw (#1_a_1) -- ++(0.25,-0.25) edge[densely dotted] +(0.5,-0.5);
							\draw (#1_a_{m-2}) -- ++(-0.25,-0.25) edge[densely dotted] +(-0.5,-0.5);
							\ifthenelse{\equal{#1}{k}}{
								\node[point] (#1_z) at (5.5,-5) [label=right:$z_k$] {};
							}{
								\node[point] (#1_z) at (5.5,-5) [label=left:$z_{#1}$] {};
							}
							\draw (#1_w_0) -- (#1_z) -- (#1_w_m);
						\end{scope}
					}

					\begin{scope}[xshift=-230,yshift=100]
						\sawedtreeblock{1}
					\end{scope}

					\begin{scope}[xshift=-81,yshift=100]
						\sawedtreeblock{2}
					\end{scope}

					\begin{scope}[xshift=100,yshift=100]
						\sawedtreeblock{k}
					\end{scope}

					\begin{scope}[xscale=0.3, yscale=0.4]
						\draw[every node/.style=point] (1_w_m) 
							-- ++(1,1) node (p1_a_0) {}
							-- ++(1,-1) node (p1_w_0) {}
							-- ++(1,1) node (p1_a_1) {}
							-- ++(1,-1) node (p1_w_1)  {}
							-- ++(0.25,0.25)
							edge[densely dotted] +(0.5,0.5)
							++(2.25,0.25) coordinate (p1c)
							++(2.25,0.25)
							edge[densely dotted] +(-0.5,-0.5)
							-- ++(0.25,0.25) node {}
							-- ++(1,-1) node (p1_w_{m-1})  {}
							-- ++(1,1) node (p1_a_{m-1}) {}
							-- (2_w_0);
						\node at (p1c) {$\cdots$};
						\foreach \x in {0,1,{{m-1}}}
						{
							\draw[every node/.style=point] (p1_w_\x)
								-- ++(0,-1) node {}
								-- ++(0,-1) node {}
								-- ++(0,-0.25)
								edge[densely dotted] +(0,-0.5)
								++(0,-1) coordinate (p1_m_\x)
								++(0,-1.5)
								edge[densely dotted] +(0,0.5)
								-- ++(0,-0.25) node (p1_b_\x) {};
							\node at (p1_m_\x) {$\vdots$};
						}
						\draw (p1c) 
							++(0,-2) node {$\cdots$}
							++(0,-2) node {$\cdots$};
					\end{scope}

					\begin{scope}[xscale=0.3, yscale=0.4]
						\draw[every node/.style=point] (2_w_m) 
							-- ++(1,1) node {}
							-- ++(1,-1) node (p2_w_0) {}
							-- ++(0.25,0.25)
							edge[densely dotted] +(0.5,0.5);
						\draw[every node/.style=point] (p2_w_0)
							-- ++(0,-1) node {}
							-- ++(0,-1) node {}
							-- ++(0,-0.25)
							edge[densely dotted] +(0,-0.5)
							++(0,-1) coordinate (p2_m_0)
							++(0,-1.5)
							edge[densely dotted] +(0,0.5)
							-- ++(0,-0.25) node (p2_b_0) {};
						\node at (p2_m_0) {$\vdots$};
						\draw (p2_w_0)
							++(2.5,0.5)
							++(0,-2) node {$\cdots$}
							++(0,-2) node {$\cdots$};
					\end{scope}

					\begin{scope}[xscale=0.3, yscale=0.4]
						\draw[every node/.style=point] (k_w_0) 
							-- ++(-1,1) node {}
							-- ++(-1,-1) node (pk_w_{m-1}) {}
							-- ++(-1,1) node {}
							-- ++(-0.25,-0.25)
							edge[densely dotted] +(-0.5,-0.5);
						\draw[every node/.style=point] (pk_w_{m-1})
							-- ++(0,-1) node {}
							-- ++(0,-1) node {}
							-- ++(0,-0.25)
							edge[densely dotted] +(0,-0.5)
							++(0,-1) coordinate (pk_m_{m-1})
							++(0,-1.5)
							edge[densely dotted] +(0,0.5)
							-- ++(0,-0.25) node (pk_b_{m-1}) {};
						\node at (pk_m_{m-1}) {$\vdots$};
						\draw (pk_w_{m-1})
							++(-2.5,0.5)
							++(0,-2) node {$\cdots$}
							++(0,-2) node {$\cdots$};
					\end{scope}

					\node[point] (z) at (-2,0) [label=below:$\bot$] {};

					\draw (z)
						edge (1_z)
						edge (p1_b_0)
						edge (p1_b_1)
						edge (p1_b_{m-1})
						edge (2_z)
						edge (p2_b_0)
						edge (pk_b_{m-1})
						edge (k_z);

					\begin{scope}[xscale=0.3, yscale=0.4]
						\draw (1_w_0) ++(0,4) node (l1_w_0) {$u_1$};
						\draw (1_w_m) ++(-2.5,4) node (l1_w_m) {$v_1$};
						\draw (p1_a_0) ++(-0.5,3) node (lp1_a_0) {$w_{1,0}$};
						\draw (p1_w_0) ++(1.5,4) node (lp1_w_0) {$w_{1,1}$};
						\draw (p1_a_1) ++(3.5,3) node (lp1_a_1) {$w_{1,2}$};
						\draw (p1_a_{m-1}) ++(0,3) node (lp1_a_{m-1}) {$w_{1,m_1}$};
						\draw (2_w_0) ++(2.5,4) node (l2_w_0) {$u_2$};
						\draw (2_w_m) ++(0,4) node (l2_w_m) {$v_2$};
						\draw (k_w_0) ++(0,4) node (lk_w_0) {$u_k$};
						\draw (k_w_m) ++(0,4) node (lk_w_m) {$v_k$};
						\draw[-Latex,gray] (l1_w_0) edge (1_w_0);
						\draw[-Latex,gray] (l1_w_m) edge[bend left] (1_w_m);
						\draw[-Latex,gray] (lp1_a_0) edge (p1_a_0);
						\draw[-Latex,gray] (lp1_w_0) edge[bend right] (p1_w_0);
						\draw[-Latex,gray] (lp1_a_1) edge (p1_a_1);
						\draw[-Latex,gray] (lp1_a_{m-1}) edge (p1_a_{m-1});
						\draw[-Latex,gray] (l2_w_0) edge[bend right] (2_w_0);
						\draw[-Latex,gray] (l2_w_m) edge (2_w_m);
						\draw[-Latex,gray] (lk_w_0) edge (k_w_0);
						\draw[-Latex,gray] (lk_w_m) edge (k_w_m);
					\end{scope}
						
				\end{tikzpicture}
			\end{equation*}
			\caption{Construction of $F'$ from $G_1, \ldots, G_k$ and the paths $p_1, \ldots p_{k-1}$.}
			\label{fig:F prime construction inductive step; proof:frames of PL p-morphic images of sawed trees}
		\end{figure}
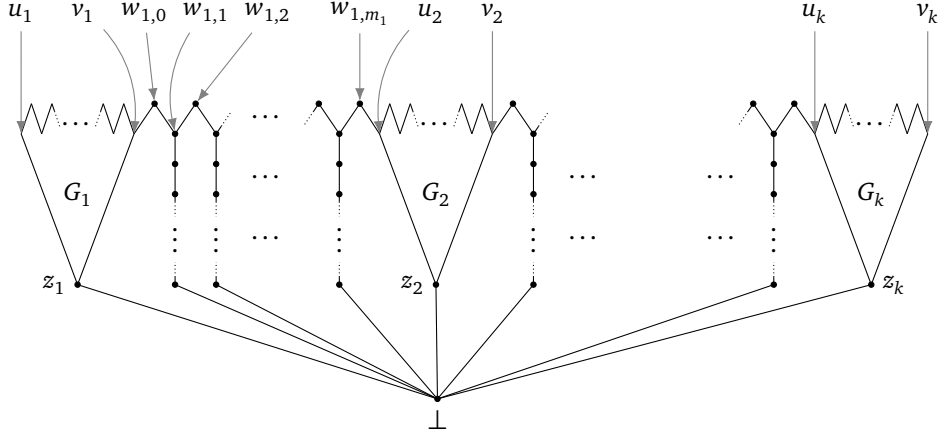

		The result is evidently a sawed tree. Finally, construct the p-morphism $f \colon F' \to F$ as follows.
		\begin{enumerate}[label=(\alph*)]
			\item Inside each sawed tree $G_i$, let $f$ act as $g_i$.
			\item For each $w_{i,j}$, let $f(w_{i,j}) \coloneqq a_{i,j}$.
			\item For each $w_{i,j}$ with $j$ odd, send the rope ladder hanging below $w_{i,j}$ to $a_{i,j}$. \qedhere
		\end{enumerate}
	\end{proof}

	\begin{corollary}\label{cor:PLn logic of sawed trees}
		$\PL_n$ is the logic of sawed trees of height at most $n$, for every $n \geq 2$.
	\end{corollary}

	\begin{proof}
		This follows from \cref{lem:sawed trees satisfy PLn} and \cref{lem:PLn frames specific connectedness property}, and the fact that $\PL_n$, like any intermediate logic, is the logic of its rooted frames.
	\end{proof}

\subsection{Convex geometric realisation}

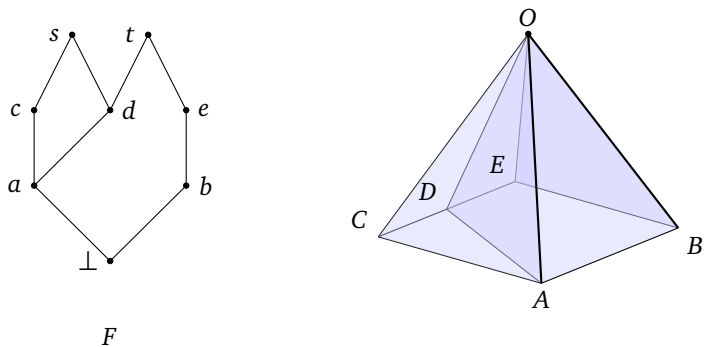
\begin{figure}[ht]
		\begin{equation*}
			\begin{tikzpicture}[scale=1]
				\tdplotsetmaincoords{110}{320}
				\begin{scope}
					\graph[poset=1] {
						o[x=1cm, label=left:$\bot$] -- {
							x0[label=left:$a$] -- {
								x[label=left:$c$] -- a[x=0.5,label=left:$s$],
								y[label=right:$d$] -- {a,b[x=0.5,label=left:$t$]}
							},
							z0[label=right:$b$] -- z[label=right:$e$] -- b
						};
					};
					\node at (1,-1) {$F$};
				\end{scope}
				\begin{scope}[xshift=180,yshift=30,scale=0.7,tdplot_main_coords]
					\coordinate (a0) at (0,0,0);
					\coordinate (a2) at (0,2,0);
					\coordinate (a4) at (0,4,0);
					\coordinate (b0) at (4,0,0);
					\coordinate (b2) at (4,2,0);
					\coordinate (b4) at (4,4,0);
					\coordinate (t) at  (2,2,4);
					\graph[use existing nodes] {
						a0 -- a2 -- a4;
						t -- {a0, a2, a4};
						b0 -- a0;
						b4 -- {a2,a4};
					};
					\fill[facefill] (t) -- (a0) -- (b0);
					\fill[facefill] (t) -- (a2) -- (b4);
					\fill[facefill] (t) -- (a4) -- (b4);
					\fill[facefill] (t) -- (b0) -- (b4);
					\graph[use existing nodes] {
						b0 -- b4;
						{[edges={thick}] {b0, b4} -- t};
					};
					\node[point] at (t) {};
					\node [above=0.05 of t] {$O$};
					\node [below=0.05 of b4] {$A$};
					\node [below right=0.01 of b0] {$B$};
					\node [above left=0.01 of a4] {$C$};
					\node [above left=0.01 of a2] {$D$};
					\node [above left=0.01 of a0] {$E$};
					\begin{scope}[tdplot_screen_coords]
					\end{scope}
				\end{scope}
			\end{tikzpicture}
		\end{equation*}
		\caption{A height-$3$ example of convex geometric realisation}
		\label{fig:convex geometric realisation height 3}
	\end{figure}

	In the second stage of the completeness proof, we provide a method of constructing a convex realisation of any sawed tree. To provide intuition for the construction, we first examine an instructive example of height $3$. Consider \cref{fig:convex geometric realisation height 3}. 

	The sawed tree $F$, depicted on the left, is realised in the pyramid $P = \mathit{OABEC}$, depicted on the right. The point $\mathit{D}$ lies midway between $\mathit{C}$ and $\mathit{E}$. An open surjective polyhedral map $f \colon P \to F$ is then defined as follows.
	\begin{itemize}
		\item The point $\mathit{O}$ is mapped to $\bot$.
		\item The remainder of the line $\mathit{OA}$ is mapped to $a$ while the remainder of $\mathit{OB}$ is mapped to $b$.
		\item The remainder of the triangle $\mathit{OAC}$ is mapped to $c$, the remainder of $\mathit{OAD}$ is mapped to $d$, and the remainder of $\mathit{OBE}$ is mapped to $e$.
		\item Finally, the remainder of the region $\mathit{OACD}$ is mapped to $s$ and the remainder of the region $\mathit{OABED}$ is mapped to $t$.
	\end{itemize}
	It is clear that such a map is polyhedral. Further the construction ensures that any open neighbourhood in $P$ is mapped to an upwards-closed subset of $F$. For instance, note that any open set intersecting $\mathit{OAD}$ must also intersect $\mathit{OACD}$ and $\mathit{OABED}$. Hence, $f \colon P \to F$ is an open polyhedral map as required.

	Notice that the two middle layers $(a,b)$ and $(c,d,e)$ of $F$ correspond to the edges $\mathit{AB}$ and $\mathit{CDE}$ of the base of the pyramid. Note further that the preimage of the tree part of $F$ --- i.e. the union of the triangles $\mathit{OAC}$, $\mathit{OAD}$ and $\mathit{OBE}$ --- has a natural triangulation. The definition of $f$ on this region then follows just as in the definition of the geometric realisation from \cref{ssec:geometric realisation}, with respect to this triangulation.

	With this intuition in mind we proceed with the proof in full generality. We make use of the following technical lemma on nerves and simplicial complexes.

	\begin{lemma}\label{lem:nerve geometric realisation criterion}
		Let $F$ be a poset and take any function $\alpha \colon F \to \R^n$. The collection:
		\begin{equation*}
			\{\Conv \alpha[X] \mid X \in \N(F)\}
		\end{equation*}
		forms a simplicial complex if and only if $\Conv\alpha[X]$ and $\Conv\alpha[Y]$ are disjoint for any disjoint $X,Y \in \N(F)$.
	\end{lemma}

	\begin{proof}
		This follows from \cite[Theorem~2]{demendez1999}, noting that the nerve $\N(F)$ is in particular an abstract simplicial complex, as defined there, with vertex set $\{\{x\} \mid x \in F\}$. 
	\end{proof}

	\begin{proof}[Proof of \cref{thm:PLn completeness for CPn}]
		The case $n=0$ is immediate. For $n=1$ note that by \cref{lem:meaning of PLn on posets}:
		\begin{equation*}
			\PL_1 = \Logic(\FPoint, \FOneFork, \FTwoFork) = \Logic(\FTwoFork)
		\end{equation*}
		Consider the convex polyhedron given by the interval $[0,1]$. We can define an open polyhedral map $f \colon [0,1] \to \FTwoFork$ by mapping $\half$ to the root, and the intervals $[0,\half)$ and $(\half,1]$ to each top node, respectively. Therefore:
		\begin{equation*}
			\Logic(\PolyCon_1) \sse \Logic([0,1]) \sse \PL_1 
		\end{equation*}
		Hence we may assume that $n \geq 2$. By \cref{cor:PLn logic of sawed trees} and \cref{lem:maps duality}, it suffices to show that every sawed tree of height $n$ can be realised geometrically in a convex polyhedron of dimension $n$. So, let $F$ be a height-$n$ sawed tree based on $(T,\prec)$. Using \cref{cor:sawed trees are planar}, let $d$ be a plane drawing of $F$ such that $d_2(x) = \height(x)$ for each $x \in F$.

		We first construct a simplicial complex corresponding to the tree part $T$ of $F$. Let $e_0, \ldots, e_n$ be the standard basis vectors of $\R^{n+1}$. Define a function $\alpha \colon T \to \R^{n+1}$ by letting, for $x \in T$:
		\begin{equation*}
			\alpha(x) \coloneqq e_{\height(x)} + d_1(x)e_{n}
		\end{equation*}
		It is helpful to consider the $n$th dimension (spanned by $e_n$) as running from left to right. Then nodes which are further to the right in the plane drawing $d$ map to points which are further to the right in $\R^{n+1}$. For each $X \in \N(T)$, let:
		\begin{equation*}
			\sig(X) \coloneqq \Conv\alpha[X]
		\end{equation*}
		Note that each element in $X$ is of a different height, so that $\alpha[X]$ is an affinely independent set of points; hence $\sig(X)$ is a simplex. Then set:
		\begin{equation*}
			\Sig \coloneqq \{\sig(X) \mid X \in \N(X)\}
		\end{equation*}
		Let us use \cref{lem:nerve geometric realisation criterion} to verify that $\Sig$ is a simplicial complex. Take disjoint $X,Y \in \N(F)$, and suppose for a contradiction that $\sig(X) \cap \sig(Y) \neq \es$. Let $X = \{x_1, \ldots, x_k\}$ and $Y = \{y_1, \ldots, y_l\}$, enumerated according to the order $<$ on $T$. Then, using barycentric coordinates inside $\sig(X)$ and $\sig(Y)$, there must be $r_1,\ldots, r_k \geq 0$ and $q_1, \ldots, q_l \geq 0$ with $\sum_{i=1}^k r_i = 1$ and $\sum_{j=1}^l q_j = 1$ such that:
		\begin{equation*}
			\sum_{i=1}^k r_i \alpha(x_i) = \sum_{j=1}^l q_j \alpha(y_j)
		\end{equation*}
		Using the definition of $\alpha$ and the fact that $e_0, \ldots, e_n$ are linearly independent, we see that:
		\begin{itemize}
			\item $r_i = 0$ if there is no $y_j$ with $\height(x_i) = \height(y_j)$,
			\item $q_j = 0$ if there is no $x_i$ with $\height(x_i) = \height(y_j)$,
			\item $r_i = q_j$ whenever $\height(x_i) = \height(y_j)$, and
			\item $\sum_{i=1}^k r_i d_1(x_i) = \sum_{j=1}^l q_j d_1(y_j)$.
		\end{itemize}
		Hence, we may assume that $k=l$ and that $\height(x_i) = \height(y_i)$ for each $i$. Now, for each $i$, since $X$ and $Y$ are disjoint, we must have $d(x_i) \neq d(y_i)$. But, since $d_2(x_i) = \height(x_i) = d_2(y_i)$, we must have either $d_1(x_i) < d_1(y_i)$ or $d_1(x_i) > d_1(y_i)$. Without loss of generality, assume that $d_1(x_1) < d_1(y_1)$. Then, since $T$ is a tree and no edges overlap in the plane drawing $d$, we must have $d_1(x_i) < d_1(y_i)$ for each $i$. Thus:
		\begin{equation*}
			\sum_{i=1}^k r_i d_1(x_i) = \sum_{i=1}^l q_i d_1(x_i) < \sum_{j=1}^l q_j d_1(y_j)
		\end{equation*}
		which is a contradiction. Therefore, $\Sig$ is a simplicial complex. As in \cref{ssec:geometric realisation}, the p-morphism $\max \colon \N(T) \to T$ gives rise to an open polyhedral map $f_T \colon \abs\Sig \to T$.

		Let us turn our attention now towards the top part of $F$. Enumerate $\Top(T)$ according to $\prec$ as $\{t_1, \ldots, t_k\}$, and let $s_1, \ldots, s_{k-1}$ be the top elements of $F$, as in the definition of a sawed tree. For each $i \leq k$, we have the $(n-1)$-simplex $\tau_i \coloneqq \sig(\ds{t_i})$. For $i \leq k-1$, let:
		\begin{equation*}
			\xi_i \coloneqq \Conv(\alpha[\Ds{s_i}]) = \Conv(\tau_{i-1} \cup \tau_i)
		\end{equation*}
		By considering the definition of $\alpha$, and noting that $\Ds{s_i}$ contains two elements which have the same height, we can see that $\Dim(\xi_i) = n$. Note also that:
		\begin{equation*}
			\xi_i \cap \xi_{i+1} = \tau_i
		\end{equation*}

		Define $P \coloneqq \bigcup_{i=1}^k \xi_i$, which will be our convex geometric realisation. By \cref{lem:dimension of union}, $P$ is an $n$-dimensional polyhedron. Furthermore, note that:
		\begin{equation*}
			P = \Conv (\tau_1 \cup \tau_k) = \Conv(\tau_1 \cup \cdots \cup \tau_k) = \Conv(P)
		\end{equation*}
		so that $P$ is a convex polyhedron and thus $P \in \PolyCon_n$. Extend the map $f_T$ to $f \colon P \to F$ by letting $x \in \xi_i \setminus (\tau_{i-1} \cup \tau_i)$ map to $s_i$. This map is clearly polyhedral. To see that it is open, take $x \in P$ and $U \sse P$ a small open neighbourhood of $x$. There are two cases. If $x \in \xi_i \setminus (\tau_{i-1} \cup \tau_i)$ for some $i$, then (as long as $U$ is small enough), $f[U] = \{s_i\}$ which is open. Otherwise, $x \in \tau_i$ for some $i$. Since $f_T$ is open, $V \coloneqq f[U \cap \abs\Sig]$ is an open subset of $T$. To see that $f[U]$ is open then, it suffices to show that whenever $s_i \in \uset^F V \cap \Top(F)$, we have $U \cap \xi_i \neq \es$. So take such an $s_i$. Since $V$ is open in $T$, we must have $t_{i-1} \in V$ or $t_i \in V$. Without loss of generality, assume the former. Hence we must have $U \cap \tau_{i-1} \neq \es$. But then since $U$ is open, it follows that also $U \cap \xi_i \neq \es$.

		Thus $f \colon P \to F$ is an open surjective polyhedral map from a convex $n$-dimensional polyhedron, as required.
	\end{proof}

\section{Conclusion}
\label{sec:conclusion}

In this article, we have provided an axiomatisation of the logic of the class of convex polyhedra. This result fits into a natural programme of investigation, initiated in \cite{tarski-polyhedra} and continued in \cite{polycompleteness}, which seeks to map out the landscape of polyhedrally complete logics. 

\begin{figure}
    \begin{equation*}
        \begin{tikzpicture}
            \begin{scope}[yshift=250]
                \node[point] (bd0) at (0,0) [label=left:{$\CPC=\BD_0$}] {};
            \end{scope}
            \begin{scope}[yshift=190]
                \node[point] (bd1) at (0,0) [label=left:$\BD_1$] {};
                \node[point] (pl1) at (7,1.5) [label=right:{$\PL_1=\LF_2$}] {};
                \draw (pl1)
                    -- ++(-1.05,-0.225) node[point] [label=below:$\LF_3$] {}
                    -- ++(-1.05,-0.225) node[point] [label=below:$\LF_4$] {}
                    -- ++(-1.05,-0.225) node[point] [label=below:$\LF_5$] {}
                    -- ++(-0.5025,-0.1125)
                    edge[densely dotted] +(-0.5025,-0.1125);
                \draw (bd1)
                    -- ++(0.25125,0.05625)
                    edge[densely dotted] +(0.25125,0.05625);
            \end{scope}
            \begin{scope}[yshift=150]
                \node[point] (bd2) at (0,0) [label=left:$\BD_2$] {};
                \node[point] (pl2) at (7,1.5) [label=right:$\PL_2$] {};
                \node[point] (flat2) at (3.5,0.75) [label=below:$\Flat_2$] {};
                \draw (bd2) 
                    edge[bend left=20] (pl2)
                    edge[bend right=20] (pl2);
                \draw (bd2)
                    -- (flat2) 
                    edge[bend left=20] (pl2)
                    edge[bend right=20] (pl2);
                \node (i2) at (-1.5,0.8) {Infinite};
                \draw[-Latex, bend left=10] (i2) edge (1.5,0.5);
            \end{scope}
            \begin{scope}[yshift=100]
                \node[point] (bd3) at (0,0) [label=left:$\BD_3$] {};
                \node[point] (pl3) at (7,1.5) [label=right:$\PL_3$] {};
                \draw (bd3) 
                    edge[bend left=20] (pl3)
                    edge[bend right=20] (pl3);
                \node (i3) at (-1.5,0.8) {Infinite};
                \draw[-Latex, bend left=10] (i3) edge (1.5,0.5);
            \end{scope}
            \begin{scope}[yshift=0]
                \node[point] (ipc) at (0,0) [label=left:$\IPC$] {};
                \node[point] (pl) at (7,1.5) [label=right:$\PL$] {};
                \node[point] (sl) at (5.025,1.125) [label=below:$\SL$] {};
                \draw (ipc) 
                    edge[bend left=20] (pl)
                    edge[bend right=20] (pl);
                \node (iw) at (-1.5,0.8) {Infinite};
                \draw[-Latex, bend left=10] (iw) edge (1.5,0.5);
            \end{scope}
            \draw (bd0)
                -- (bd1)
                -- (bd2)
                -- (bd3)
                -- ++(0,-0.5)
                edge[densely dotted] +(0,-0.5)
                ++(0,-1.25) node {$\vdots$}
                ++(3.5,0.75) node {$\vdots$};
            \draw (ipc)
                -- ++(0,0.25)
                edge[densely dotted] +(0,0.25);
            \draw (pl1)
                -- (pl2)
                -- (pl3)
                -- ++(0,-0.5)
                edge[densely dotted] +(0,-0.5)
                ++(0,-1.25) node {$\vdots$};
            \draw (pl)
                -- ++(0,0.25)
                edge[densely dotted] +(0,0.25);
            \draw (pl1) -- (bd0);
        \end{tikzpicture}
    \end{equation*}
    \caption{The currently-mapped landscape of polyhedrally complete logics. $\CPC$ is classical logic: $\IPC$ plus the principle of excluded middle.}
    \label{fig:landscape of poly-complete logics}
\end{figure}
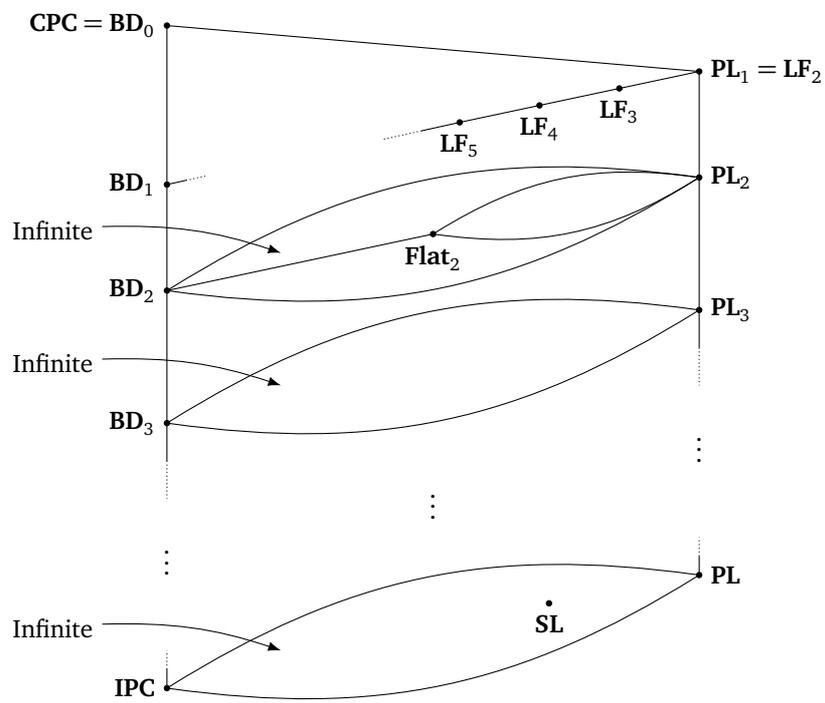

In \cite{polycompleteness} it is shown that there are infinitely many polyhedrally complete logics of each height, axiomatised by the Jankov-Fine formulas of `starlike trees'. This in particular includes Scott's logic $\SL$. Beyond these results, \cite{gabelaiatacl} investigates the lower-level structure of this landscape in more detail. First, it is shown that every height-$1$ logic is polyhedrally complete: these are $\BD_1$ plus the logic $\LF_k$ of the `$k$-fork' --- the frame consisting of a root with $k$ immediate successors --- for each $k \geq 2$. Second, turning to the height-$2$ case, the focus is on logics of `flat polygons': $2$-dimensional polyhedra which can be embedded in the plane $\R^2$. Any such logic turns out to be axiomatised by a \emph{subframe formula} (see \cite[p.~313]{chagrovzakharyaschev1997}) plus the Jankov-Fine formulas of certain trees. Moreover, there is a smallest such logic: $\Flat_2$.
\Cref{fig:landscape of poly-complete logics} charts out what is currently known about the landscape of polyhedrally complete logics, to the best of our knowledge. 

One long-term goal is the complete classification of all polyhedrally complete logics. This article presented one schema for attacking this problem: starting with a natural class of polyhedra and asking what its logic is. For this it is important to be able to find a geometric realisation of any frame of a candidate logic in the class of polyhedra under consideration. By contrast, in \cite{polycompleteness} another schema is followed. There we start from the logic side and define a class of logics with the aim that they are polyhedrally complete, making use of the Nerve Criterion for polyhedral completeness.

\vspace{3mm}

\noindent
{\bf Acknowledgement} The authors would  like to acknowledge support by the SRNSF Grant \#FR-22-6700.

	\printbibliography

\end{document}